\documentclass{article}
\usepackage[margin=1in]{geometry}
\usepackage{authblk}
\usepackage{caption}
\usepackage[style=nature,backend=biber,autocite=superscript]{biblatex}

\usepackage[dvipsnames]{xcolor}
\DeclareBibliographyCategory{highlight}

\AtEveryBibitem{
  \ifcategory{highlight}{\color{blue}}{}
}

\definecolor{NMblue}{HTML}{025E8D}
\usepackage[
  colorlinks   = true,
  linkcolor    = NMblue,
  citecolor    = NMblue,
  urlcolor     = NMblue
]{hyperref}
\usepackage{amsfonts,amssymb,amsmath,mathrsfs}
\usepackage{graphicx}
\usepackage{soul}
\usepackage{bm}
\usepackage{enumerate,algorithmicx,algorithm}
\usepackage{changepage}
\usepackage{comment}

\usepackage{tabularx}
\usepackage{hhline}
\usepackage{diagbox}

\usepackage{array,booktabs}
\usepackage{siunitx}
\graphicspath{{./figures}}

\newcommand{\be}{\begin{equation}}
\newcommand{\ee}{\end{equation}}
\newcommand{\ba}{\begin{aligned}} 
\newcommand{\ea}{\end{aligned}}
\newcommand{\bea}{\begin{eqnarray}} 
\newcommand{\eea}{\end{eqnarray}}
\newcommand{\mbf}[1]{{\mathbf #1}}

\newcommand{\eps}{\varepsilon}
\newcommand{\qqquad}{\qquad\qquad}

\newcommand{\R}{\mathbb{R}}

\newcommand{\Z}{\mathbb{Z}}

\newcommand{\rr}{{\mathbf{r}}}
\newcommand{\x}{\mbf{x}}

\newcommand{\kk}{{\bm{\xi}}}
\newcommand{\ki}{\mbf{k}}
\newcommand{\n}{\mbf{n}}

\newcommand{\F}{\mbf{F}}
\newcommand{\LL}{{\mathcal L}}
\newcommand{\Sp}{{\mathcal S}}
\newcommand{\hS}{{\hat {\mathcal S}}}
\newcommand{\tp}{\tilde{\phi}}
\newcommand{\bl}{{\bm{\ell}}}
\newcommand{\al}{\alpha}
\newcommand{\rc}{{r_c}}
\newcommand{\bigO}{{\mathcal O}}

\DeclareMathOperator{\erf}{erf}
\DeclareMathOperator{\erfc}{erfc}

\newcommand{\tcb}[1]{\textcolor{black}{#1}}

\newcommand{\os}{\hspace{1ex}}

\begin{document}
\title{Accelerating Molecular Dynamics Simulations using Fast Ewald Summation with Prolates}

\author[1,2]{Jiuyang Liang}
\author[1]{Libin Lu}
\author[1,*]{Alex Barnett}
\author[1,3,*]{Leslie Greengard}
\author[1,*]{Shidong Jiang}

\affil[1]{Center for Computational Mathematics, Flatiron Institute, Simons Foundation, New York, New York 10010, USA}
\affil[2]{School of Mathematical Sciences, Shanghai Jiao Tong University, Shanghai, 200240, China}
\affil[3]{Courant Institute of Mathematical Sciences, New York University, New York, New York 10012, USA}

\date{\today}

\maketitle

\begingroup
\renewcommand{\thefootnote}{\fnsymbol{footnote}}
\footnotetext[1]{Correspondence: abarnett@flatironinstitute.org, lgreengard@flatironinstitute.org, sjiang@flatironinstitute.org}
\endgroup

\begin{abstract}
The evaluation of long-range Coulomb interactions is a significant cost in molecular
dynamics (MD), even when using Particle Mesh Ewald (PME) or 
Particle-Particle-Particle-Mesh (PPPM) methods, which rely on Ewald splitting and 
the fast Fourier transform to achieve near-linear scaling. 
We introduce ESP---Ewald summation with prolate spheroidal wave functions (PSWFs)---which
leads to a more efficient Fourier representation
and a reduction in the required grid size, global communication, and particle-grid 
operations, without loss of accuracy.
We have integrated the ESP method into two widely-used open-source
MD packages, LAMMPS and GROMACS, enabling rapid comparison and adoption.
Relative to PME/PPPM baselines 
at error tolerances $10^{-3}$ to $10^{-4}$,
ESP gives roughly a $3$-fold acceleration of electrostatic interactions,
and a $2.5$-fold speed-up in the MD simulation when
using about $10^3$ compute cores.
At high accuracy ($10^{-5}$), these increase to $10$-fold
for the far-field electrostatics and $5$-fold for MD simulation.
Furthermore, we show that the accelerated codes have improved strong scaling with core count, and
validate them in realistic long-time biological and material simulations.
ESP thus offers a practical, drop-in path to reduce the time-to-solution and energy footprint of MD workflows.
\end{abstract}


Molecular dynamics (MD) is an indispensable tool in both classical and 
quantum physics, supporting investigations of material properties, chemical
reactions, protein structure, and drug design in both academia and industry~\autocite{karplus1990molecular,karplus2002molecular,souza2021martini,fu2022accurate,ives2024restoring}. 
MD has proven particularly powerful for the quantitative characterization of physicochemical
and biomolecular systems, complementing experimental data, guiding experimental design, and predicting properties of systems that are costly or challenging to probe experimentally.
While MD simulations trace atomic interactions and trajectories on the 
femtosecond time scale, the phenomena of interest typically unfold over
microseconds, milliseconds or even longer. 
As a result, typical simulations involve
$10^9$ to $10^{12}$ time steps, requiring that each time step be executed 
in well under one
millisecond if results are to be obtained within days or weeks.
Achieving this performance makes parallel computing and the efficient use 
of large CPU/GPU clusters essential.

At the atomic scale, MD simulations involve both
short-range forces, due to bond vibrations, torsional rotation and
Lennard--Jones interactions, as well as 
long-range Coulomb forces. While the computational cost of computing
all short-range forces is of the order $O(N)$ for $N$ atoms, the
Coulomb interactions require $O(N^2)$ work using direct computation.
In the literature, there are two main families of fast algorithms that 
overcome this cost:
tree-based methods~\autocite{barnes1986nature,greengard1994fast}, such as
the Fast Multipole Method (FMM)~\autocite{greengard1987jcp,greengard1988} 
and its descendants, which can achieve $O(N)$ complexity and
handle highly inhomogeneous particle
distributions; and Fourier-based methods that rely on the FFT (the
fast Fourier transform)~\autocite{fft}.
While the asymptotic scaling of the latter methods is of the order 
$O(N \log N)$, the prefactor is smaller than for the tree-based schemes 
for relatively homogeneous systems. As this is the case for solid and 
liquid phase atomic systems,
Fourier-based methods are generally considered the methods of choice.
The first such scheme was proposed by Ewald~\autocite{ewald1921ap} before the advent of the FFT,
achieving an $O(N^{3/2})$ complexity 
through ``kernel splitting'' and the Poisson summation formula.
In the last several
decades, a wide variety of faster schemes have been introduced that differ
in technical detail but follow the same general principle. 
While we do not seek to review the literature here,
these $O(N \log N)$ schemes include PME (Particle Mesh Ewald)~\autocite{darden1993jcp,deserno1998jcp,essmann1995jcp,shan2005jcp},
Particle-Particle-Particle-Mesh (PPPM or P$^3$M)~\autocite{hockney1988,ballenegger2008jcp}, and their many variants, which are typically optimized to achieve
3 to 5 significant digits in the force evaluation.
PME and PPPM are the default modules for computing Coulomb interactions
 in leading open-source MD codes such as
LAMMPS~\autocite{thompson2021lammps},
 GROMACS~\autocite{berendsen1995gromacs}, and NAMD~\autocite{schulten2020}.
In contexts where greater precision is required, the spectral Ewald method achieves higher-order accuracy~\autocite{specewald,shamshirgar2021jcp}. Spectral Ewald methods have also been incorporated into GROMACS~\autocite{SEGromacs}.

A drawback of all Ewald/Fourier methods, however, is that on 
modern parallel architectures, they require global communication 
(all-to-all) of $O(N)$ data for each three-dimensional (3D) FFT.
This is because the latter involves global data
transpose operations between the stacked 1D FFTs.
Thus, a significant fraction of the total cost is spent on Coulomb
interactions (this varies from $25\%$ in GROMACS default settings,
to $>80\%$ in our LAMMPS experiments),
and this fraction grows with the system size, or as higher accuracy is demanded.
Recent attempts to mitigate the communication cost include the 
Random Batch Ewald
method~\autocite{jin2021sisc},
which replaces the global gather-scatter with randomized sub-samples. 
While effective in a statistical sense, its accuracy is guaranteed only on 
average and may fall short for applications that demand strict error control.

In Ewald's original formulation~\autocite{ewald1921ap}, the 3D Coulomb kernel is 
decomposed into two parts:
\begin{equation}
 \frac{1}{r} \;=\; \frac{ \erf(r/\lambda)}{r} + \frac{ \erfc(r/\lambda)}{r}, 
\label{ewalderf}
\end{equation}
where $\lambda$ is a length scale (typically of order the interatomic spacing in MD),
and
\[  \erf(r) := \frac{2}{\sqrt{\pi}} \int_0^r e^{-t^2} \, dt,\qquad 
\erfc(r) := 1 - \erf(r). \]
The first part---involving an antiderivative of a Gaussian---is long-range
but smooth, and well-approximated using only Fourier frequencies up to some bandlimit.
The second part, by contrast, is singular but short-range, so that spatial truncation to a radius of a few $\lambda$ incurs very small error.
Here, we propose a kernel-splitting strategy
using the zeroth-order prolate
spheroidal wave function (PSWF)~\autocite{slepian1961bstj,landau1961bstj,slepian1978bstj,slepian1983sirev,osipov2013} in place of the Gaussian.
It has been known in
harmonic/Fourier analysis and
signal processing since the 1960s that PSWFs possess an optimal energy-concentration property:
specifically, among all square-normalized functions supported in a spatial
interval, a PSWF carries the maximum possible energy within a given frequency bandlimit.
The relevance for computation is that this splitting permits a much shorter
(in fact optimally short) Fourier series to represent the long-range, periodic, portion of the electrostatic field, without any increase in truncation radius or error. At high accuracy, keeping the cutoff radius fixed yields nearly a twofold reduction in grid size per dimension; in 3D this corresponds to an $\approx 8\times$ savings in the FFT stage alone.
This results in faster FFTs and lower communication overhead in high-performance computing (HPC) settings.

A secondary benefit of PSWFs is that they can minimize the
cost of {\em particle-grid} operations, namely (a)
spreading charges from particle locations to nearby points on a uniform grid and
(b) interpolating potentials or forces from the grid back to the particles themselves.
Unlike PME and PPPM, which typically use B-spline kernels, a
Cartesian product of PSWFs
yields better accuracy at every MD step with fewer operations, and more reliable
error control in potential and force evaluation. Using PSWFs for both kernel splitting and particle-grid operations further cuts the cost of the long-range Coulomb interactions. By contrast, PME/PPPM reach the same accuracy only with many more grid points per dimension unless heavy upsampling is used; in practice, their spreading/interpolation order is typically set to four or five to limit particle-grid cost, which in turn forces much longer FFTs. As a result, even at modest accuracy (three to four digits), PSWFs reduce the Fourier-transform length by roughly $6$--$14\times$. In highly optimized codes such as GROMACS, where dynamic load balancing often assigns $\sim\! 75\%$ of the work to short-range interactions, the cutoff radius can be reduced to rebalance short- and long-range work, yielding an overall $2$--$3\times$ speedup.

We call our method, which uses
PSWFs to replace both the split in Eq.~\eqref{ewalderf} and the
spreading/interpolation kernels, {\em Ewald summation with prolates} (ESP).
To facilitate a fair benchmark comparison with the classical fast Ewald methods, we
made minimal modifications to both the
state-of-the-art LAMMPS and GROMACS codes, replacing only those function modules
relevant to electrostatics while leaving all other components untouched.
For instance, within LAMMPS, ESP results in an overall
three- to seven-fold acceleration over PPPM with recommended
settings, and better strong scaling with respect to
number of cores, while achieving the same four-digit force accuracy.
Within GROMACS, ESP allows a 195 ns/day simulation speed
for 1M atoms on 960 cores with a 2 fs timestep,
at an error tolerance of $2\times 10^{-4}$, a $2.5\times$ speed-up over the already highly-optimized
GROMACS code.
The acceleration---which is mostly due to the fact that the PSWF enables coarser Fourier grids and thus faster FFTs---
depends strongly on the required error tolerance:
it grows at higher accuracy, where the PPPM/PME methods require much more upsampling of the Fourier grid.
Even at the lower 3-digit accuracy typical for biological applications, we show a three-fold acceleration of the
Coulomb calculation within GROMACS.
Given that MD simulations currently account for more
than 20\% of the core-hours on the top-500
supercomputers~\autocite{nersc2020,antypas2013,bottaro2024,pall2014,pall2022},
such improvements translate into significant reductions in computational time, energy consumption, and operational cost. ESP thus has the potential for a broad impact on large-scale molecular modeling.

\begin{table}
  {\scriptsize
  \setlength{\tabcolsep}{4.pt}
  \begin{tabular}{|llllllllll|}
      \hline
  system & $N$ atoms & $\Delta$ (error)& code & cores & PME or PPPM & ESP  & speed-up & Fig. & notes
  \\
&&&&& $h$(nm), $\rc$(nm), $P$& $h$(nm), $\rc$(nm), $P$ &&&
  \\
  \hline
  \multicolumn{10}{|l|}{\rule{0pt}{3ex} \hspace{-2ex} {\bf Coulomb only, lower accuracy:}}\\
  water &   375--12M &  $10^{-3}$ & G & 1 &    0.13, \os 1.0, 4    &       0.14, \os 0.5, 4 &         3.2$\times^*$  &  1a & \\
  water &  375--12M &  $10^{-4}$ & G & 1  &    0.083, 1.0, 4   &       0.1, \os\os 0.5, 5 &          3.3$\times^*$  &  1b & \\
  membrane  &  82k    &  $4\times 10^{-4}$ &  G & 12--384
    &  0.12, \os 1.0, 4  &  0.2, \os\os 0.7,  5 &         2.6$\times^*$  &  1c  &  MEM in \cite{kutzner2025scaling} \\
  ribosome  & 2M & $10^{-3}$ & G & 96--3k
  &  0.135, 1.0, 4   &      0.16, \os 0.6, 4  &        3$\times^*$  &    1d  &  RIB in \cite{kutzner2025scaling} \\
trunc. octa. & 62k & $4\times10^{-4}$ & G & 1--96 & 0.12,\os\os 1.2, 4 & 0.12,\os\os 0.5, 5 & 4.1$\times^*$ & S1a & CONA in \cite{musleh2024analysis}\\
rhomb. dodec. & 1.2M & $2\times10^{-4}$ & G & 12--768 & 0.09,\os\os 1.0, 4 & 0.16,\os\os 0.7, 5 & 4$\times^*$ & S1b & RNA in \cite{posani2025ensemble}\\
  \hline
    \multicolumn{10}{|l|}{\rule{0pt}{3ex} \hspace{-2ex} {\bf Large MD, higher accuracy:}}\\
water  & 3.6M, 106M & $10^{-4}$ & L & 12--9k &   0.10, \os 0.9, 5    &     0.26, \os 0.9, 5  &   3--7$\times$ &  2ab & \\
water  & 4M, 192M  & $10^{-5}$  & G & 12--9k &   0.046, 0.9, 5    &    0.19, \os 0.9, 8  &  5--6$\times$  &  S2ab & \\
\hline
    \multicolumn{10}{|l|}{\rule{0pt}{3ex} \hspace{-2ex} {\bf Long-time MD, higher accuracy:}}\\
lysozyme & 1M  & $2\times 10^{-4}$ & G & 960 &   0.12, \os 0.9, 5   & 0.25, \os 0.9, 5   &   2.5$\times$  & 3b & \\
LiTFSI & 1M  &  $10^{-5}$        & G   & 960 &   0.049, 0.9, 5  &   0.194, 0.9, 8    &  5$\times$ & 4b &  \\
\hline
  \end{tabular}
  }
  \caption{Summary of all timing performance comparisons reported in this paper,
    grouped into the three main types of study listed in Results.
    PME, particle mesh Ewald; PPPM, particle-particle-particle-mesh; ESP, Ewald summation with prolates; MEM, membrane benchmark; RIB, ribosome benchmark; CONA, Concanavalin A benchmark; RNA, rhombic-dodecahedral RNA benchmark.
    $\Delta$ is the targeted relative force error as defined in \eqref{Deltadef};
    PME/PPPM (default) or ESP parameters are adjusted to achieve this error.
    The codes are either LAMMPS or GROMACS, indicated by L or G.
    $h$ is the grid spacing (``Fourier spacing''), and $r_c$ is the Coulomb cutoff radius.
    In the full molecular dynamics simulations, the Coulomb and Lennard--Jones cutoff radii are identical, whereas in Coulomb-only benchmarks the Lennard--Jones interactions are disabled. 
     $P$ is the number of grid points per dimension (``order'') in particle-grid interactions. 
    Speed-ups are
    for the full molecular dynamics simulation, apart from
    those with an asterisk which are for the Coulomb part only. In figure references,
    ``S'' denotes supplementary figure.
    Error and other validation experiments are shown in Supplementary Figs.~3--7.
}    
    \label{tbl:expts}
\end{table}

\section*{Results}

We assess the speed of the proposed method using the set of
timing benchmarks summarized in Table~\ref{tbl:expts}.
Each row corresponds to an experiment, or set of experiments, on one type of molecular system,
and lists the chief parameter settings.
In each case,
to ensure a fair comparison,
the parameters were hand-tuned to achieve the stated relative force error with both the PME/PPPM
and the ESP method.
Specifically for the existing PME/PPPM methods, the default recommended spreading
and interpolation kernel width (order) $P$ and ``Coulomb'' radius $r_c$ are used, while
the Fourier grid spacing $h$
is decreased until the target force relative accuracy
$\Delta$ (see Methods) is achieved; such parameters are labeled as ``Default'' throughout this paper.

Our benchmarks group into \tcb{three types} of study (matching the three blocks of Table~\ref{tbl:expts}):

\begin{enumerate}
\item {\bf Coulomb-only timings at lower accuracy.}
  We start by measuring the Coulomb calculation speed-up factor achieved by integrating our ESP method into the GROMACS package; see Fig.~\ref{fig:time}.
  We target electrostatic force relative errors $\Delta$ in the range $10^{-3}$ to $10^{-4}$, as
  typically required in biological applications.
  The single-core performance is measured using bulk water systems of up to 12M atoms.
  The multi-core strong scaling is measured using a membrane and a ribosome system from the recent
  benchmark paper of Kutzner et al.\ \autocite{kutzner2025scaling},
  together with a truncated-octahedral Concanavalin A (CONA) system \autocite{musleh2024analysis} and a rhombic-dodecahedral RNA system \autocite{posani2025ensemble}. For these benchmarks, we adopt the default
    PME parameters and target errors $\Delta$ from the original studies.
  We find an acceleration of the Coulomb evaluation by roughly a factor of three.
  
\item
  {\bf Large MD simulations at higher accuracy and core count.}
  Here large-scale water systems with between 3.5M and 192M atoms are
  simulated for short times.
  We measure the MD simulation
  speed-up factor
  that results by inserting our ESP method into both LAMMPS and GROMACS packages.
  Within LAMMPS, targeting an error $\Delta=10^{-4}$,
  we study the CPU time reduction achieved by parallelizing
  over a growing number of compute cores; see Fig.~\ref{fig:LAMMPS_time}(a,b).
  We find that ESP has better strong scaling, and that
  beyond $10^3$ cores ESP results in
  an overall factor of $5$--$7$ acceleration of the MD simulation relative to the default PPPM method.
  For GROMACS, at the higher-accuracy target error $\Delta=10^{-5}$
  typically demanded in material applications~\autocite{thompson2021lammps,bagchi2020surface,antila2018dielectric},
  we find a $5$--$6$ speed-up factor.
  Details of the LAMMPS results are given below in this section, while GROMACS results are summarized in Supplementary Fig.~2.

\item
  {\bf Long-time organic MD simulations at higher accuracy.}
  Finally we present long-time
  simulations of two medium-scale (around 1M atom) organic systems, using GROMACS:
  i) lysozyme proteins in solution,
  and ii) a microphase-separated high-concentration lithium aqueous electrolyte
  popular in Li-ion batteries.
  Across both systems, our approach achieves speedups of $2.5$ and $5$ relative to GROMACS,
  at error levels of $\Delta=2\times10^{-4}$ and $\Delta=10^{-5}$, respectively.
  We also validate the simulations by showing that the change in method
  has no detectable effect on a variety of spatiotemporal and thermodynamic quantities.
\end{enumerate}

Referring the reader to the Methods section for more details on the set-up, we now
break down the results of the three types of study.

\subsection*{Performance for Coulomb calculations at lower accuracy}

Our first study focuses on the single-core CPU performance, and the parallel CPU scaling, of our ESP method for the Coulomb calculation alone, when compared to the Ewald-based PME method supplied by the software package GROMACS.
For single-core tests, we use the SPC/E bulk water system~\autocite{berendsen1987missing}
with the number of atoms $N$ between 375 and 12,288,000.
We used ``default'' PME parameters (see row 1 and 2 of Table~\ref{tbl:expts}, and Methods)
to achieve relative force errors of either $\Delta=10^{-3}$ or $10^{-4}$, covering the typical range of
tolerances in biological simulations~\autocite{darden1993jcp, berendsen1995gromacs, schulten2020, case2005amber}.
For the scaling test, we use a membrane system (81,743 atoms) which was also used to study water permeation through an aquaporin tetramer that is embedded in a lipid bilayer~\autocite{de2001water}, and a ribosome system (2,037,834 atoms), which is used for the study of the function of these protein factories in our cells~\autocite{bock2013energy}.
These two systems were also part of a recent GROMACS benchmark \autocite{kutzner2025scaling},
from which we used the default PME parameters as in row 3 and 4 of Table~\ref{tbl:expts},
and match their target force errors of $\Delta = 4 \times 10^{-4}$ and $10^{-3}$, respectively. These first four benchmark systems all use orthorhombic boxes (cubic for water, and rectangular for membrane and ribosome). To test non-orthorhombic periodic boxes, we also benchmark two systems: truncated-octahedral CONA (62,045 atoms) and rhombic-dodecahedral RNA (1,232,355 atoms). Following the default PME setups and target errors from existing benchmarks~\autocite{musleh2024analysis,posani2025ensemble}---as detailed in rows 5 and 6 of Table~\ref{tbl:expts}---we achieved accuracies of $\Delta=4\times10^{-4}$ and $2\times 10^{-4}$, respectively.

\begin{figure}[!htbp]
\centering
\includegraphics[width=0.9\linewidth]{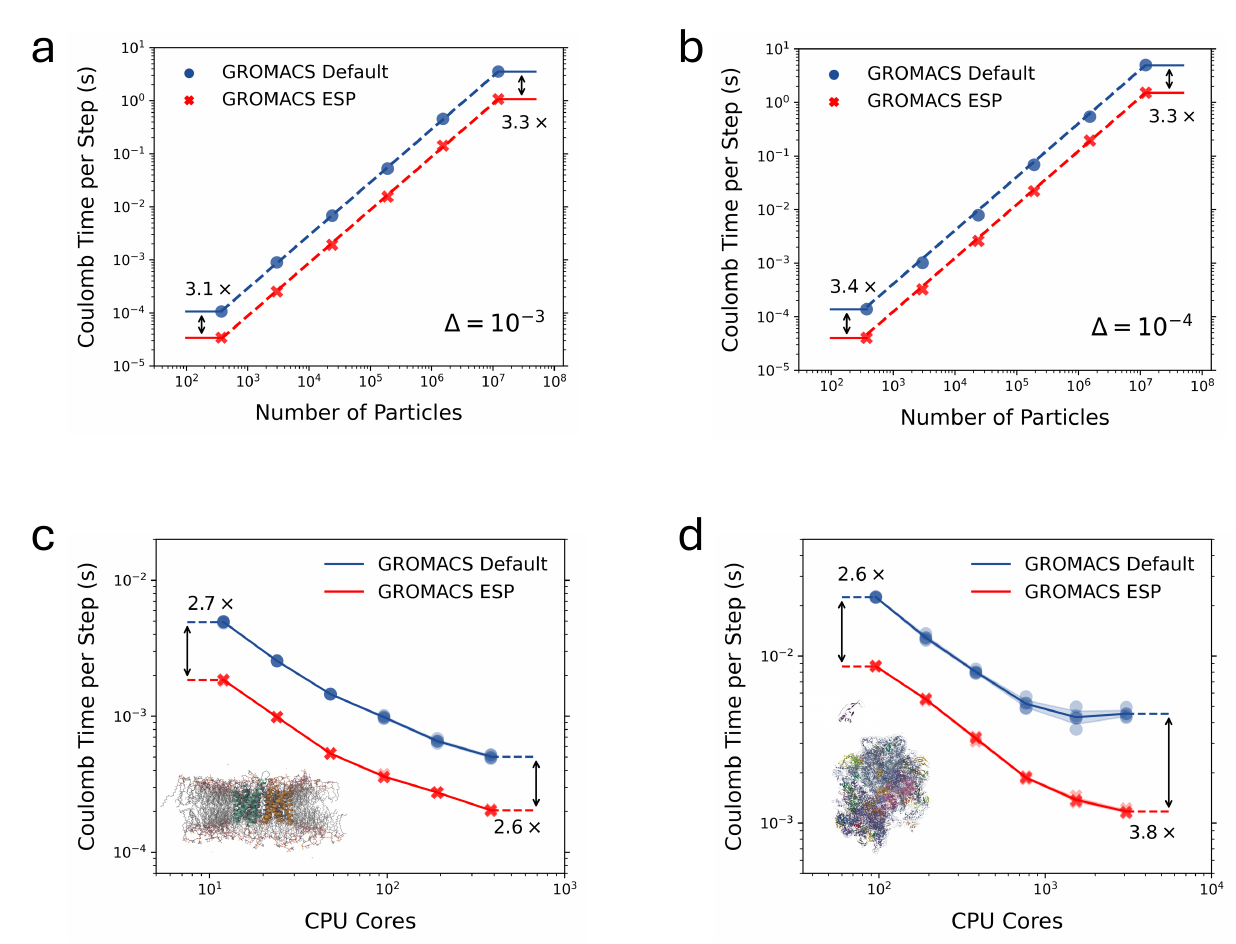}
\caption{\sf 
Performance comparison of Ewald summation with prolates (ESP) and particle mesh Ewald (PME) for Coulomb calculations in representative biomolecular benchmark systems. ({\bf{a}}, {\bf b}) Bulk water systems benchmarked on a single core. The average computational time of Coulomb interactions per simulation step, measured over 30-minute runs, is shown for two error tolerances: {\bf a}, $\Delta = 10^{-3}$; and {\bf b}, $\Delta = 10^{-4}$, plotted against the number of particles. Results were obtained using the ESP method implemented in GROMACS and compared with the default PME algorithm. Blue circles represent PME with default settings; red `\texttimes{}' symbols indicate ESP. Dashed lines represent linear fits to the data. 
({\bf{c}}, {\bf d}) Performance comparison of ESP with PME for membrane (81,743 atoms) and ribosome (2,037,824 atoms) systems. Insets show snapshots of the corresponding systems. The average Coulomb computation time per step is plotted as a function of the number of central processing unit (CPU) cores, with each point averaged over five runs of $30$ minutes. The PME reference data in Panels {\bf c} and {\bf d} follow setups reported in the Methods section, with estimated error tolerances of $\Delta=4\times10^{-4}$ and $10^{-3}$, respectively. Source data are provided as a Source Data file.}
\label{fig:time}
\end{figure}

\begin{figure}[!htbp]
\centering
\includegraphics[width=0.98\linewidth]{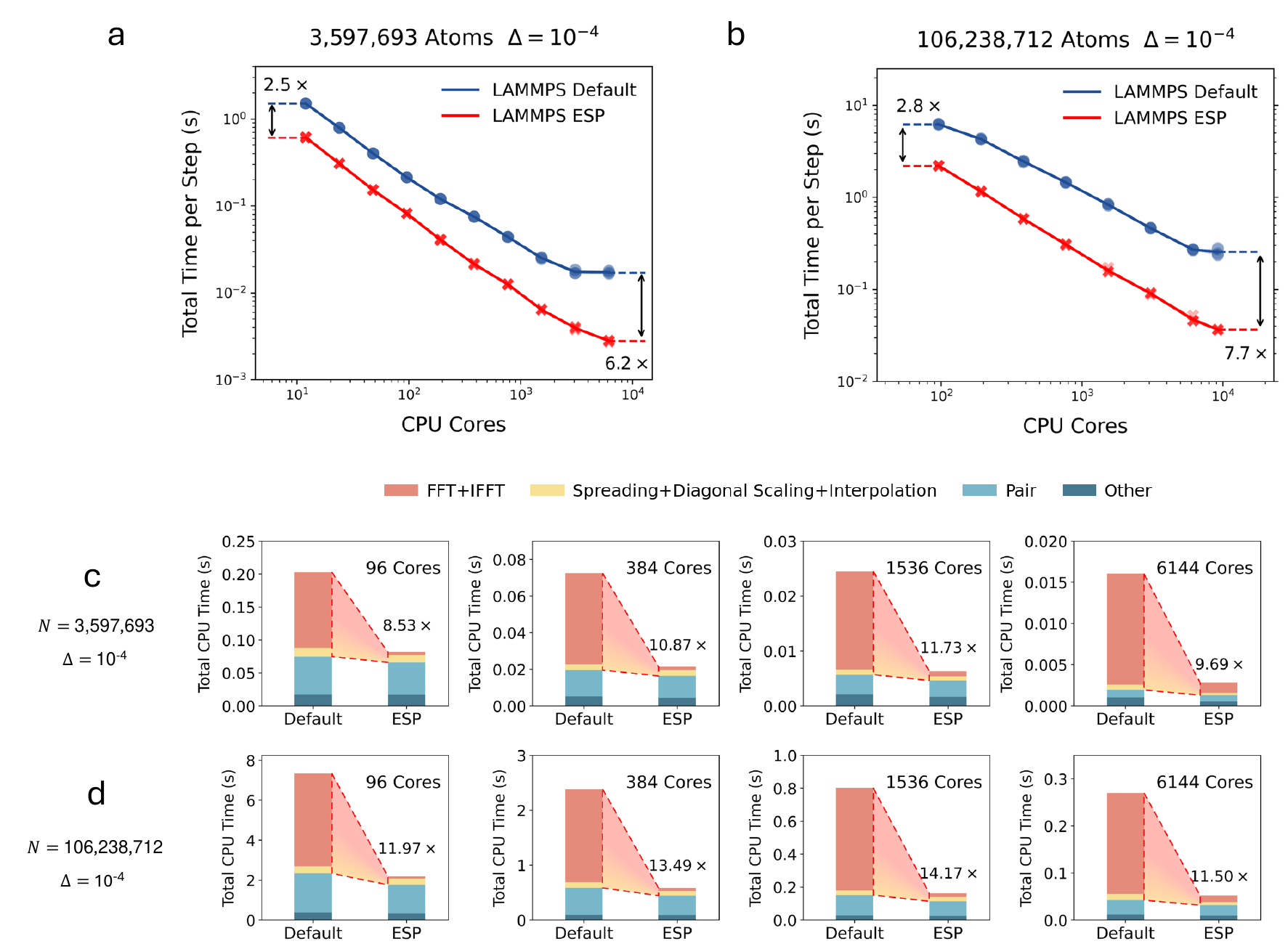}
\caption{\sf Performance comparison of Ewald summation with prolates (ESP) and particle-particle-particle-mesh (PPPM) for large bulk water molecular dynamics simulations in LAMMPS. Total simulation time per step, averaged over five runs of $30$ minutes each, is shown for systems with 3,597,693 atoms ({\bf a}) and 106,238,712 atoms ({\bf b}) as a function of the number of central processing unit (CPU) cores. Data were generated using ESP implemented within LAMMPS, compared against the native PPPM option, with an error tolerance of $\Delta = 10^{-4}$. 
Blue circles show PPPM with default parameters, and red '\texttimes{}' marks show ESP. Light-colored markers show the results of five repeated runs. Solid lines indicate the mean across the five runs, and shaded bands indicate the $95\%$ confidence intervals of the mean. In panels {\bf c} and {\bf d}, the simulation time per step is broken down into four components: fast Fourier transform (FFT) and inverse fast Fourier transform (IFFT) operations (red); spreading, diagonal scaling, and interpolation (yellow); short-range pairwise interactions (light blue); and all other simulation tasks like thermostats and constraints (dark blue). Panels {\bf c}-{\bf d} correspond to the systems shown in {\bf{a}}-{\bf b}, respectively. Speedups in the long-range interaction are annotated directly within each panel. Source data are provided as a Source Data file.}
\label{fig:LAMMPS_time}
\end{figure}

In all six cases we isolate the computational time for Coulomb interactions from a $30$-minute run and use only the second half of the run to compute the average Coulomb time per step.
The Coulomb calculation time per step, including both the local and spectral components (see Methods), is shown in Fig.~\ref{fig:time}(a,b). Compared with the default parameter settings, ESP achieves a three-fold speed-up.
With a single core, this speed-up is independent of $N$; both methods scale linearly with $N$.
To assess parallel scaling performance, we report the Coulomb CPU time per step
for the membrane and ribosome systems as the number of CPU cores increases to 3,000.
As shown in Fig.~\ref{fig:time}(c,d) and Supplementary Fig.~1, the ESP speedup increases with core count, achieving maximum accelerations of $3.8$-fold and $5.1$-fold for orthorhombic and triclinic systems, respectively.
This indicates that ESP has improved strong scaling over PME, a finding that we believe is
explained by the significant reduction in Fourier grid size.

The speedup factor is explained by two different parameter changes, as follows.
The Coulomb radius $r_c$ is defined as the truncation
distance of the short-range part (the first term in Eq.~\eqref{ewalderf} or its PSWF variant).
The $r_c$ values of $0.5$ to $0.7$ nm that we are able to use with ESP
(see top four rows of Table~\ref{tbl:expts})
are much smaller than the default $r_c=1.0$ nm used with PME (as taken from \cite[Tbl.~1]{kutzner2025scaling}).
Our smaller $r_c$ allows the local part of the Coulomb calculation, which theoretically scales as $\bigO(r_c^3)$,
to be much cheaper.
At the same time, the grid spacing (or ``Fourier spacing'', labeled $h$ in Table~\ref{tbl:expts})
is larger with ESP than was possible with PME, making the number of Fourier grid points smaller,
hence the FFTs faster, while still matching the required force error.
These two parameter improvements ($r_c$ and $h$) combine to explain the roughly three-fold electrostatics acceleration.
Note that in ESP the particle-grid ``order'' parameter was set to $P=4$ (equal to the choice
used by \cite{kutzner2025scaling} as indicated by the TPR files linked from that paper),
or, at the slightly higher accuracies, $P=5$, in order to enable a larger $h$.

In contrast to the above, in
the remaining benchmarks, which are full MD simulations, we no longer hand-tune $r_c$ below its default value.
This is for the technical reason that $r_c$ typically also sets the truncation radius for
van der Waals (Lennard--Jones) interactions in the GROMACS or LAMMPS codes, and that $r_c < 0.9$ nm would result in
unacceptably inaccurate simulations.

\subsection*{Performance in simulations of large-scale systems at higher accuracy}
  
Next, we compare the parallel scaling and performance of our
ESP method against the Gaussian-splitting PPPM and PME methods supplied by
the software packages LAMMPS and GROMACS respectively. For this, we use the
SPC/E bulk water system~\autocite{berendsen1987missing}.
For comparison with LAMMPS at the moderate accuracy $\Delta = 10^{-4}$, we considered systems with
$N=$ 3,597,693
and 106,238,712 atoms.
To facilitate the study of repeated runs of large systems over a wide range of numbers of cores, a small number (2000) of time-steps was used.
The total time per simulation step---including short- and long-range nonbonded forces, bonded
interactions, neighbor list updates, constraints, and thermostats---is reported in
Fig.~\ref{fig:LAMMPS_time}(a,b). When compared to the default PPPM parameter settings, ESP achieves a
$3$--$7$ fold increase in MD simulation speed.
The factor increases with the core count, indicating superior strong scaling,
becoming a factor of $7$ at many thousands of cores.

In GROMACS, we now target a high force accuracy, $\Delta=10^{-5}$, common in materials science
and other applications that demand high precision, such as electrolyte solutions, ionizable biomolecules, and interfacial environments \autocite{thompson2021lammps,bagchi2020surface,antila2018dielectric}.
We chose large systems with $N=$ $4,214,784$
and $192,000,000$.
The total time
per simulation step, as shown in Supplementary Fig.~2(a,b), demonstrates that ESP consistently achieves a factor of $5$--$6$ speedup,
although this does not improve with core count.
Supplementary Fig.~3(a,b) confirms that ESP achieves the same or
better electrostatic accuracy than PPPM and PME across all tested conditions.  The finding that the speed-up from using ESP varies amongst software packages may be attributed to different levels of code optimization. To ensure a fair comparison at the algorithmic level, we have preserved the existing parallel optimization architecture of each package.

The large observed acceleration factor using ESP is explained by
its reduction in the time for the long-range Coulomb task.
This is clear in the breakdowns of the total CPU time per MD time-step for the four water systems
shown for LAMMPS in Fig.~\ref{fig:LAMMPS_time}(c-d), and for GROMACS in
Supplementary Fig.~2(c-d).
In each plot, the time spent in the forward and inverse FFT (including communication
costs) is indicated in red.
The combined cost of spreading, diagonal scaling, and interpolation is indicated in 
yellow.
These first two contributions together comprise the long-range Coulomb force
calculation.
The cost of short-range pairwise interactions, including the short-range real-space component
of the Coulomb forces, is indicated in light blue.
All other costs are indicated in dark
blue. (See Methods for definitions of these components.)
These plots show that long-range Coulomb calculation is an order of magnitude faster with ESP.

To understand this long-range acceleration,
recall that for both PPPM (LAMMPS) and PME (GROMACS) we set the spreading order to the ``default''
$P=5$, then tune the grid spacings $h$ to reach relative force errors of $\Delta = 10^{-4}$ or $10^{-5}$.
This results in the parameters in rows 5 and 6 of Table~\ref{tbl:expts}.
ESP is able to achieve the same target error with much larger grid spacings than PME/PPPM:
$h$ is 2.6 times larger at $\Delta=10^{-4}$, and 4 times larger at $\Delta=10^{-5}$.
The total number of Fourier grid points scales like $1/h^3$, so is reduced by factors of about 18 and 64,
respectively, explaining the dramatic reductions in FFT and communication costs.

Note that at $\Delta=10^{-4}$ ESP uses precisely the same order $P=5$ (and of course $r_c$) as PPPM, so that the
ability to use the much larger $h=0.26$ nm without loss of accuracy is entirely due to the change from Gaussian to PSWF based kernels.
At the highest accuracy $\Delta=10^{-5}$,
we found it advantageous with ESP to increase the order to $P=8$;
we have also explored increasing $P$ within PME, and, while it allows $h$ to be increased
somewhat, it does not close the gap with ESP
(for example, in PME the high $P=12$ allows $\Delta=10^{-5}$ to be reached with $h=0.11$ nm,
but this is still much smaller than the $h=0.19$ nm possible with ESP at $P=8$).

\subsection*{Lysozyme proteins in solution}

\begin{figure}[!htbp]
\centering
\includegraphics[width=0.98\linewidth]{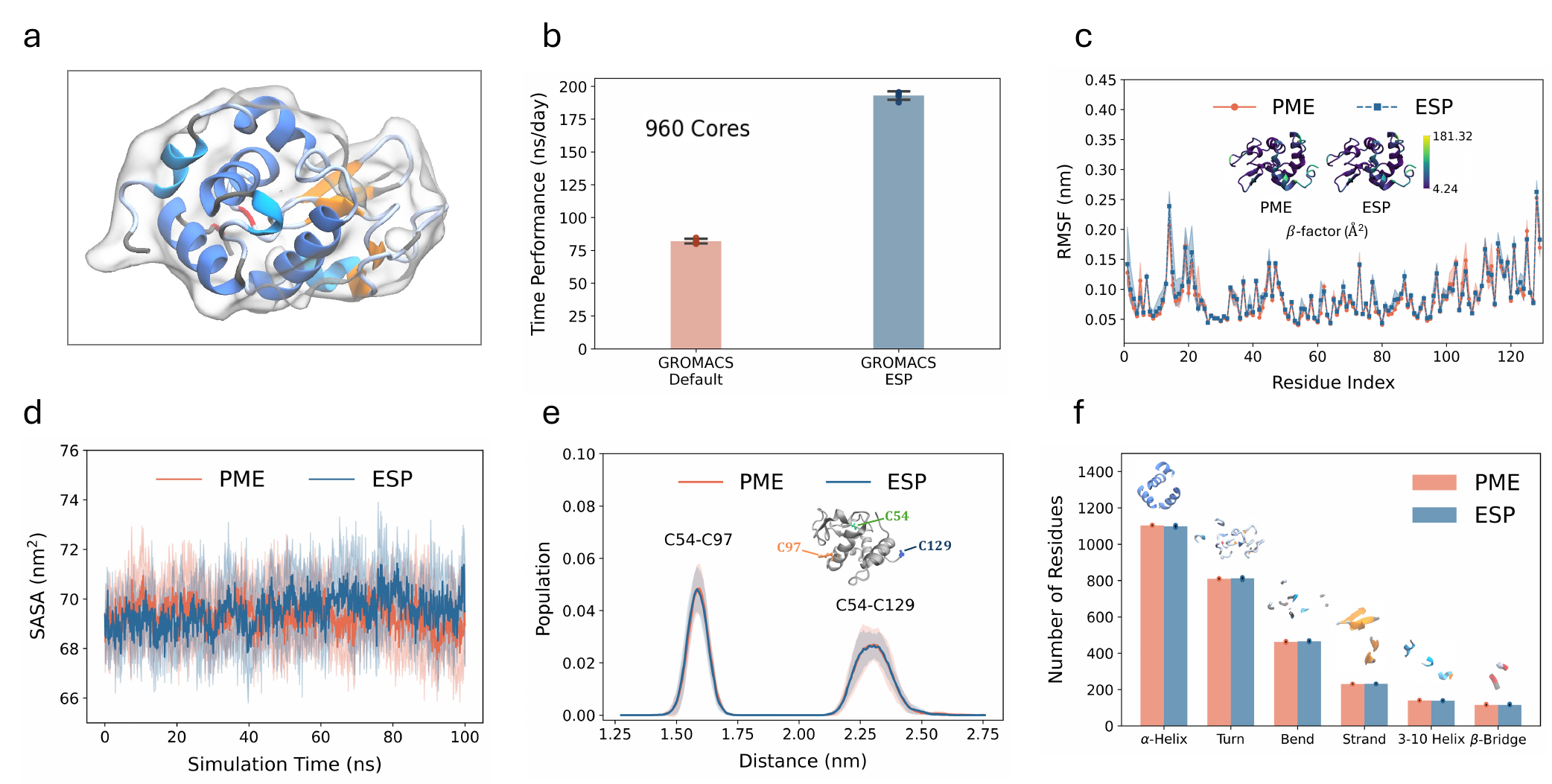}
\caption{\sf Comparison of particle mesh Ewald (PME)- and Ewald summation with prolates (ESP)-based GROMACS simulations for a large lysozyme solution benchmark. The system comprises 1,036,152 atoms, 
including $27$ duplicated lysozyme proteins. {\bf a}, Simulation snapshot of the local 
environment around one protein, with coloring based on secondary structure. 
{\bf b}, Performance comparison (see row 7 of Table~\ref{tbl:expts} for parameters). Each bar shows the mean over $5$ independent runs, and error bars indicate standard deviation.
{\bf c}, Root mean square fluctuation (RMSF) across residues; 
the inset shows two lysozyme proteins from the PME- and ESP-based simulations at the same simulation time, 
colored by residue-specific temperature ($\beta$)-factors. 
{\bf d}, Solvent-accessible surface area (SASA) averaged over all proteins during the simulation. 
Light-colored shaded areas indicate the standard deviation averaged over five independent 
runs. 
{\bf e}, Distribution of characteristic inter-domain distances between residue pairs 
C54--C97 and C54--C129. Shaded regions indicate standard deviation estimated from five runs.
{\bf f}, The total number of residues belonging to each type of secondary structure. Source data are provided as a Source Data file.}
\label{fig:protein}
\end{figure}

To benchmark the accuracy and efficiency of ESP-based MD in long-time biological
simulations, we consider a lysozyme protein system in aqueous solution---a widely used
model for testing simulation and experimental protocols. The system contains
1,036,152 atoms, including $27$ replicated lysozyme proteins; one is shown in
Fig.~\ref{fig:protein}a.
Simulations were performed for $100$ ns using $960$ CPU
cores, with ``default'' parameters $P=5$ and Fourier spacing $h=0.12$~nm for the PME. The practical error tolerance is about $2\times10^{-4}$.
Parameters are summarized in row 7 of Table~\ref{tbl:expts}.
The electrostatic errors are
validated in Supplementary Fig.~6a.
At the same error level, ESP (with the same $P$ but over twice the Fourier spacing $h=0.25$ nm)
reaches a simulation speed of approximately
$195$ ns/day, representing a $2.4$-fold speedup compared to the default GROMACS/PME methods
(see Fig.~\ref{fig:protein}b).

Biomolecular dynamics were characterized using root mean square fluctuation (RMSF) and the temperature $\beta$-factor (Fig.~\ref{fig:protein}c) across the residues of each lysozyme. Surface-accessible surface area (SASA) was analyzed to
capture characteristic structural features (Fig.~\ref{fig:protein}d). 
We also compute the root mean square deviation (RMSD) and radius of gyration (RG) of the proteins, the Gibbs free energy difference, and the Ramachandran distributions of different residues over time, with results shown in Supplementary Figs.~4 and~5.
The definitions of RMSF, SASA, RMSD, and RG are provided in Methods. 
To examine functional phase space sampling, we analyzed the distances between residue pairs C54--C97 and C54--C129, which characterize the open--close motion of the lysozyme's two 
domains involved in enzymatic activity~\autocite{faber1990mutant}
(Fig.~\ref{fig:protein}e). In addition, the number of secondary structures of each type was computed to assess protein stability (Fig.~\ref{fig:protein}f).
Using nearly $14$ times fewer FFT grid points compared to the default GROMACS setting,
GROMACS with ESP faithfully reproduced all these
structural, dynamic, and functional features observed in native GROMACS simulations.
Further details on simulation protocols and data analysis are provided in Methods.

\subsection*{Concentrated LiTFSI ionic liquid at high accuracy}

\begin{figure}[!ht]
\centering
\includegraphics[width=0.98\linewidth]{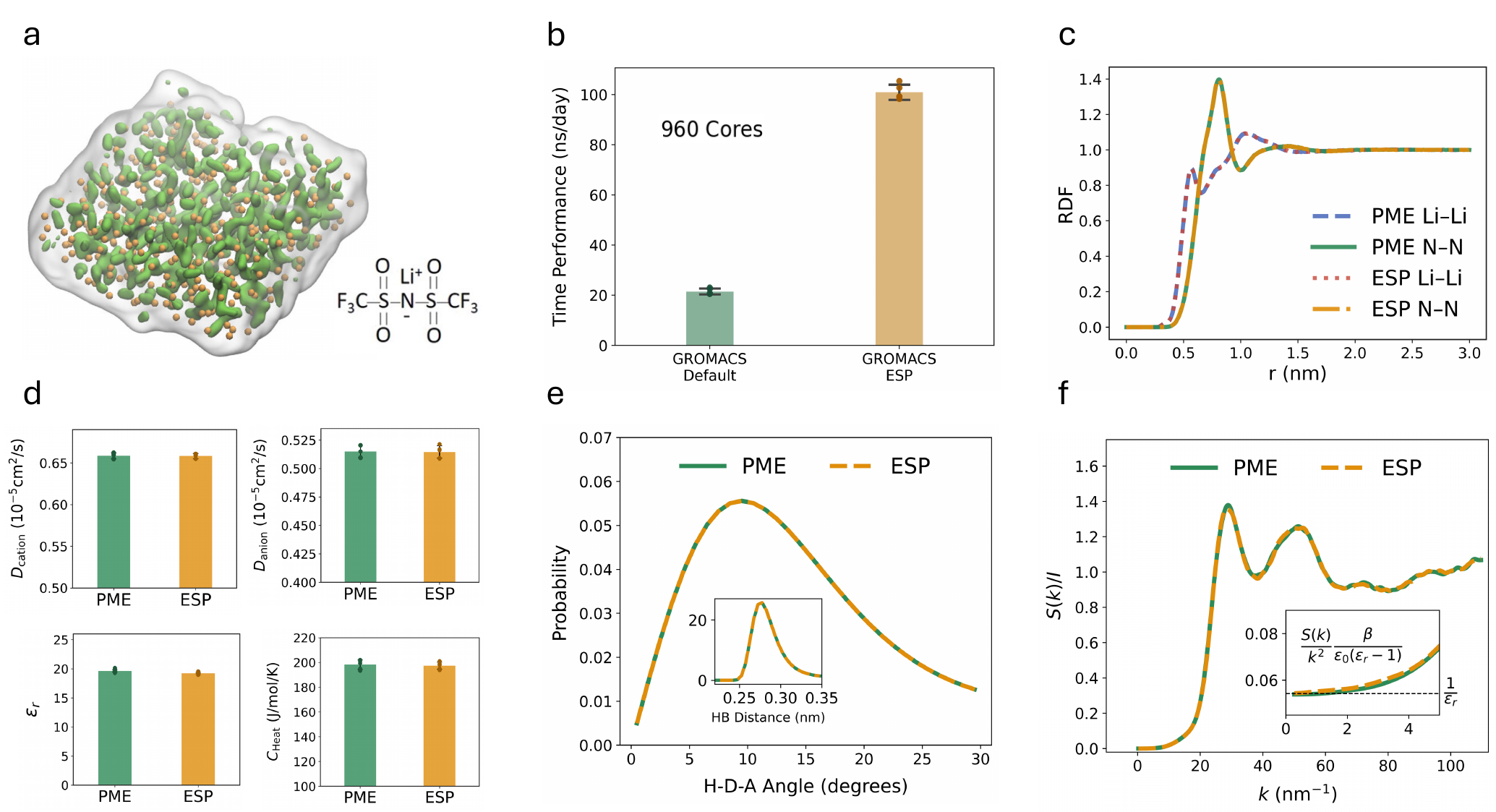}
\caption{Comparison of simulation results obtained using GROMACS with particle mesh Ewald (PME) and Ewald summation with prolates (ESP) on a concentrated Li-ion aqueous electrolyte, LiTFSI (5 M/L, $1,011,392$ atoms, using $960$ cores). {\bf a}, Simulation snapshot of the system. Cations and anions are shown in yellow and green, respectively; the solvent is represented as a white quicksurf surface.
{\bf b}, Performance comparison; see row 8 of Table~\ref{tbl:expts} for parameters. Each bar shows the mean over $5$ independent runs, and error bars indicate standard deviation.
{\bf c}, Radial distribution functions (RDF) of nitrogen--nitrogen pairs among anions and lithium--lithium pairs.
{\bf d}, Diffusion constants ($D$) of cations and anions, the dielectric constant ($\varepsilon_r$) of the aqueous electrolyte, and the system's heat capacity ($C_{\text{Heat}}$). 
{\bf e}, Probability distribution of the hydrogen--donor--acceptor (H--D--A) angle, characterizing the geometric structure of hydrogen bonds; the inset shows the distribution of hydrogen bond lengths.
{\bf f}, Charge structure factor $S(k)$ over the ionic strength $I=V^{-1}\sum_{i=1}^{N} q_i^2$ in Fourier space along the $k_z$ direction, where $\mathbf{k} = (0,0,k_z)$ and $k = |\mathbf{k}|$; the inset shows $S(k)\beta/(k^2\varepsilon_0(\varepsilon_r-1))$, with $\beta = 1/k_{\text{B}}T$ and $k_{\text{B}}$ the Boltzmann constant. The results show excellent agreement with the symmetry-preserving mean-field theory~\autocite{hu2022symmetry}, which predicts that its zero-frequency limit equals $1/\varepsilon_r$ (as calculated in {\bf d}).
Source data are provided as a Source Data file.}
\label{fig::LiTFSI}
\end{figure}

To benchmark the performance and accuracy of ESP in realistic material systems, we study the concentrated Li-ion aqueous electrolyte, lithium bis(trifluoromethanesulfonyl)imide (LiTFSI). It can be used as a Li-ion source in electrolytes for Li-ion batteries as a safer alternative to lithium hexafluorophosphate~\autocite{suo2015water}. The system has a concentration of 5~M/L, comprising $1,011,392$ atoms; see Fig.~\ref{fig::LiTFSI}a. Simulations were run for $100$~ns using 960 CPU cores, with an error tolerance of $10^{-5}$ (validated in Supplementary Fig.~6b). This level of accuracy is sufficient even for systems with stringent force-field requirements.
  The parameters are given in row 8 of Table~\ref{tbl:expts}. ESP achieves a slightly smaller error than
  the default PME, using a Fourier grid spacing $h$ that is about 4 times larger, and thus FFTs that
  are about 64 times smaller.

\tcb{ESP achieved a simulation speed of approximately $100$~ns/day, delivering a fivefold speedup over the default GROMACS setups for this accuracy (Fig.~\ref{fig::LiTFSI}b). Accuracy was evaluated based on several metrics: the radial distribution functions (RDFs) of nitrogen-nitrogen pairs among anions and lithium-lithium pairs, diffusion constants of cations and anions, the dielectric constant of the aqueous electrolyte, the system's heat capacity, the probability distribution of the hydrogen--donor--acceptor (HDA) angle, and the charge structure factor along with its zero-frequency limit. Across all metrics, ESP reproduces the structural, dynamical, thermodynamical, and dielectric correlation properties of this microscopically inhomogeneous system with accuracy similar to native GROMACS. This illustrates the large efficiency advantage of using ESP when high accuracy is needed.}

\section*{Discussion}
We have developed {\em Ewald summation with prolates} (ESP),
a fast Ewald scheme to accelerate the calculation
of long-range Coulomb interactions in molecular dynamics simulations. Instead
of the conventional Gaussian that was proposed by Ewald~\autocite{ewald1921ap} more than 100 years ago and still in use in all existing software packages, we use a radial kernel
derived from the first prolate spheroidal wave function of order
zero (PSWF) to split the kernel into a Fourier space component and a real-space
component.
This allows a larger Fourier grid spacing for the same accuracy,
hence reduces the required number of FFT grid points by a large factor.
Furthermore, we reduce the cost of the spreading and interpolation stages by replacing the B-splines introduced in 
SPME~\autocite{essmann1995jcp} (Smooth PME) thirty years ago with Cartesian-product PSWFs.
Since the PSWF is optimal in terms of frequency localization, mathematically there cannot
be any significant kernel-based improvements beyond what we propose here.
Implementing these within the popular LAMMPS and GROMACS codes involved quite localized software changes
(mostly generating sets of coefficients for use by the efficient existing
piecewise-polynomial kernel and short-range Coulomb function evaluators).
ESP can also be readily incorporated into other MD packages.

The resulting acceleration of MD simulations depends quite heavily on the required force accuracy.
At high (5-digit) accuracies, we observe speed-up factors of 5--6 within GROMACS,
whereas at moderate (four-digit) accuracies this becomes around 2.5 in GROMACS, but up to 7 in LAMMPS.
At low (3-digit) accuracies, one would expect less dramatic MD speed-up factors, because the
long-range Coulomb is less dominant.
At 3-digit accuracy we instead focused on benchmarking the acceleration of the Coulomb
evaluation alone, showing a 3-fold speed-up and improved strong scaling within GROMACS.
As well as increasing the grid spacing, the latter exploited reducing the Coulomb cutoff to as low as $r_c=0.5$ nm,
which makes the {\em short-range direct Coulomb} component much cheaper.
In order to exploit the latter in an accurate MD simulation, one would need also to
implement the van der Waals ($1/r^6$ part of the Lennard--Jones potential) within ESP.
This is because van der Waals is currently cut off at $r_c$ in the packages, so would cause unacceptable errors for
$r_c<0.9$ nm.
An ESP implementation of Lennard--Jones is in progress, and will allow us to report realistic
3-digit accuracy MD acceleration factors in the future.

While these benchmarks were performed on CPU clusters to establish a baseline against state-of-the-art implementations~\autocite{kutzner2025scaling}, the mathematical advantages of ESP are fundamentally architecture-agnostic. The reduction in Fourier grid requirements directly translates to decreased memory bandwidth pressure and reduced communication overhead---two of the primary bottlenecks in GPU-accelerated MD simulations. By optimizing the cutoff radius $r_c$, ESP provides a flexible mechanism to re-balance workloads between the real-space (particle-particle) and Fourier-space (grid-based) tasks, which is critical for maximizing throughput on heterogeneous hardware. Preliminary tests on GPU architectures already indicate a roughly two-fold speedup over standard PME, and further optimization of ESP-specific kernels for massively parallel environments is currently underway.

Our scheme can be readily applied to other boundary conditions---including quasi-1D and quasi-2D systems and those with dielectric interfaces~\autocite{mazars2011long,liang2020harmonic}---as well as to various ensembles~\autocite{frenkel2001}.
It is equally straightforward to extend it to treat other interaction kernels, and the development of optimal PSWF-based splitting for Lennard--Jones and other non-Coulombic potentials is currently under investigation.
Recent work has leveraged the PME method to train machine-learning models for more accurate long-range force-field evaluation~\autocite{gong2025predictive,ji2025mMLIPLong}; our ESP method can be seamlessly integrated into those workflows to reduce training cost.
Additionally, the ESP framework is extensible to multipolar force fields and QM/MM simulations~\autocite{giese2016jctc,nam2005jctc}, where its superior Fourier-space representation can mitigate the high computational costs of higher-order derivatives and quantum-classical coupling.
Combining our approach with the random-batch Ewald method can also 
enable more efficient stochastic sampling without increasing the sample size. 
Finally, because the ESP method significantly reduces the computational and communication overhead of Coulomb interactions while maintaining full numerical accuracy at each time step, it can be effectively combined with coarse-grained models and enhanced sampling techniques to further accelerate and scale up cross-scale simulations~\autocite{souza2021martini,yang2019enhanced,bussi2020using}.

\section*{Methods}
\subsection*{A general Ewald split for electrostatic interactions}

Consider a charged system of $N$ particles located at $\rr_{i}$, $i=1,\cdots,N$, with
charge strengths $q_{i}$, in the periodic domain $\Omega=[-L/2,L/2]^3$;
we present the cubic case for simplicity.
A core calculation in MD simulations is the evaluation of 
the electrostatic potentials $u(\rr_i)$ at these locations, 
excluding the self-contribution:
\be
u_i:= u(\rr_i)=
\sum_{\substack{j\in\{1,\dots,N\}, \; \mathbf{n}\in\Z^3 \\ (j,\mathbf{n})\neq(i,\mbf{0})}}\frac{q_j}{4\pi
  |\rr_{i}-\rr_j-\mathbf{n}L|},
\qquad i=1,\dots,N,
\label{task}
\ee
the total electrostatic energy $E=\frac{1}{2}\sum_{i=1}^{N}q_iu_i$,
and the electrostatic force acting on the $i$th particle
$\mathbf{F}(\mathbf{r}_i)=-\nabla_{\mathbf{r}_i}E$.
Here we have nondimensionalized
so that physical constants such as $\epsilon_0$ do not appear.
For periodic systems, it is well known that 
the system must satisfy the charge neutrality condition $\sum_{j=1}^{N}q_j=0$.
The potential $u$ is only determined up to an arbitrary additive constant~\autocite{de1980simulation,hu2014jctc}.

Rather than the traditional Ewald split \eqref{ewalderf}, the ESP method
considers a general radially-symmetric split of the Coulomb kernel into
spectral (smooth long-range) plus local (short-range) parts:
\be
\frac{1}{4\pi r}=\Sp(r) + \LL(r)
:=
\frac{\int_0^{r/\rc}\chi_\al(x)dx}{4\pi r} + 
\frac{1-\int_0^{r/\rc}\chi_\al(x)dx}{4\pi r},
\label{kernelsplitting}
\ee
where $\rc$ will define a cutoff (truncation) radius,
and $\alpha>0$ is a ``shape parameter'' depending on the desired precision $\eps$.
Here $\chi_\al$ is a smooth, nonnegative, even function with normalization
$\int_0^1 \chi_\al(x)dx = 1 - \bigO(\eps)$,
and $\eps$-support of $[-1,1]$, meaning $\chi_\al(1) \approx \eps$
and $\chi_\al(x)$ decays rapidly to zero as $x\rightarrow \infty$.
We will drop the subscript $\alpha$ when it is clear from context.
The short-range component $\LL(r)$
may be truncated at $r=\rc$ to precision $\eps$.
This allows the corresponding ``local'' part of the potential
\be
u_i^l = \sum_{\substack{j\in\{1,\dots,N\}, \; \n\in\Z^3 \\ (j,\n)\neq(i,\mbf{0})}} \LL (|\rr_{i}-\rr_j-L\n|)
\, q_j
\label{localpart}
\ee
to be computed directly in $O(Ns)$ cost, where $s$ is the average number of particles
within distance $\rc$ of a particle.

The long-range component $\Sp$, on the other hand, is smooth,
so it has a rapidly decaying 3D Fourier transform
\be
\hS(\kk) := \int_{\R^3} \Sp(\x) e^{i\kk\cdot\x} d\x
\;= \;
\frac{\hat{\chi}(\rc |\kk|)}{|\kk|^2}, \quad \kk \in \R^3,
\label{hS}
\ee
where $\hat{\chi}$ is the 1D Fourier transform of $\chi$ regarded as
a function on $\mathbb{R}$.
This useful connection of 3D to 1D Fourier transforms
follows by using spherical coordinates and radial integration by
parts (see App.~A.3 of the DMK framework~\autocite{jiang2025dmk}). From the above,
the long-range (``spectral'') part of
the periodic potential is 
\be
u_i^s \;:= \!\!\!\!  \sum_{j\in\{1,\dots,N\}, \; \n\in\Z^3} \!\!\!\!
    \Sp (|\rr_{i}-\rr_j-\mathbf{n}L|)\,q_j
        \;=
    \sum_{\ki\in \Z^3, \ki\neq\mbf{0}} \!
    e^{-i\frac{2\pi\ki}{L}\cdot\rr_i} \frac{\hS(2\pi\ki/L)}{L^3}
\sum_{j=1}^N e^{i\frac{2\pi\ki}{L}\cdot\rr_j} q_j,
\label{spectralpart}
\ee
where the wave vectors $\kk$ in \eqref{hS} are related to integer indices in \eqref{spectralpart}
by $\kk = 2\pi \ki / L$. 
In practice, the spectral grid in Fourier space is truncated to the index set
$I := \{-n_f/2,-n_f/2+1,\dots,n_f/2-1\}^3 \subset \Z^3$, i.e.,
a cubical set of $N_f:= n_f^3$ total Fourier modes.
The even number of modes in each dimension, $n_f$, must be chosen large enough that
$\hS(2\pi\ki/L)$ has decayed to $\eps$ everywhere for $\ki$ outside of $I$.
Recalling \eqref{hS},
this decay is controlled by the smoothness of $\chi$, the choice of $\rc$,
and (more weakly) the Coulomb factor $1/|\kk|^2$.
Note that the self-interaction term is not excluded in the spatial sum because it merely
changes the constant in the potential. The truncation error of the spectral part of the potential is defined as
\be
E_T =
    \left|\sum_{\ki\in \Z^3\setminus I}
    e^{-i\frac{2\pi\ki}{L}\cdot\rr_i} \frac{\hS(2\pi\ki/L)}{L^3}
\sum_{j=1}^N e^{i\frac{2\pi\ki}{L}\cdot\rr_j} q_j\right|.
\label{truncationerror}
\ee
\subsection*{Fast algorithm for the spectral part}

The form of \eqref{spectralpart} suggests one possible method to compute it, via a three-step
procedure: i) the rightmost exponential sum is computed using a so-called
nonuniform FFT (NUFFT) of type 1,
then ii) diagonal scaling by $\hS$ is done, finally iii) the left-most exponential sum
is computed via a NUFFT of type 2.
For details about the NUFFT see~\autocite{nufft2,nufft3,nufft6}.
While this method achieves an overall $\bigO(N + N_f \log N_f)$ complexity,
is convenient, and is used by some researchers~\autocite{Hedman2006ENUF},
it sacrifices efficiency because internally the NUFFTs
require {\em upsampling} to an FFT grid of length greater than $n_f$,
in order to guarantee accurate NUFFT outputs at the highest Fourier frequencies.
Yet the rapid decay of $\hS$ means that such accuracy is not needed in those modes,
and upsampling may be avoided, as follows.

The more efficient fast algorithm for the spectral part may be interpreted as
``reassembling'' components of the NUFFT,
namely its spreading, interpolation, FFT, and inverse FFT operations.
This generalizes the popular PME and PPPM methods in a more mathematical framework.
We present the case without upsampling.
Suppose that $\phi(\x)$ is a compact window
function used in spreading nonequispaced data $(q_j,\rr_j)$ to its
neighboring $P\times P \times P$ physical grid points, and in interpolating the potential from
the same grid points back to the particles; we call $\phi$ the ``SI kernel.''
Let $\tp(\x) := \sum_{\n\in\Z^3} \phi(\x + L\n)$ be its 3D periodic extension.
Let $\{ p_\ki\}_{ \ki \in I }$ be a set of diagonal scaling coefficients, to be determined
shortly.
Then, the set of spectral potentials $u_i^s$ can be evaluated by the following
five-step
procedure:
\bea
\ba
b_\bl &= \sum_{j=1}^N \tp(\rr_j - h\bl) q_j
,\quad \bl \in I
&\qqquad \mbox{(spread)}
\label{spread}
\\
\hat{b}_\ki &= \sum_{\ell\in I} e^\frac{2\pi i \ki \cdot \bl}{n_f} b_\bl
,\quad \ki\in I
&\qqquad \mbox{(FFT)}
\\
\hat{c}_\ki &= p_\ki \hat{b}_\ki, \quad \ki\in I
&\qqquad \mbox{(diagonal scaling)}
\\
c_{\bl'} &= \sum_{\ki\in I} e^\frac{-2\pi i \ki \cdot \bl'}{n_f} \hat{c}_\ki
,\quad \bl' \in I
&\qqquad \mbox{(IFFT)}
\\
\tilde u^s_i &= \sum_{\bl'\in I} \tp(\rr_i - h\bl') c_{\bl'}
,\quad i \in \{1,\dots,N\}
&\qqquad \mbox{(interpolate)}.
\label{interp}
\ea
\eea
Here $h=L/n_f$ is the grid spacing in physical space
(known as ``Fourier spacing'' in GROMACS), with
values listed in Table~\ref{tbl:expts}.
In the third step a good choice for $p_\ki$
(sometimes called the ``influence function''~\autocite{hockney1988,deserno1998mesh}) is
\be
p_\ki = \frac{\hS(2\pi\ki/L)}{L^3\, |\hat\phi(2\pi\ki/L)|^2}
=
\frac{\hat\chi(2\pi\rc|\ki|/L)}{2\pi L\, |\hat\phi(2\pi\ki/L)|^2\, |\ki|^2},
\quad \ki \in I,
\label{pk}
\ee
which simultaneously carries out convolution with the smoothed Coulomb kernel
and cancels the effect of the convolutions with the SI kernel $\phi$ in steps
1 and 5.
The computational cost of the five-step
method is $\bigO(N P^3 + N_f\log N_f)$:
the first term corresponds to spreading and interpolation, and the second to the FFT pair.
To make the spreading/interpolation efficient in practice, the SI kernel
is usually chosen as the Cartesian product of a 1D window function,
in turn approximated with piecewise polynomials on subintervals of length
$h$.

The forces $\F_i = -q_i \nabla u(\rr_i)$ are gradients of the potentials,
and thus inherit a split into local plus spectral parts.
Their spectral part can be computed by the five-step method above,
either via multiplying $\hat{c}_\ki$ by $2\pi i \ki/L$ in
step 3 (the ``$i\ki$'' method), or by taking analytic derivatives of
$\tp$ in step 5 (the ``AD'' method).
LAMMPS PME with default settings uses the $i\ki$ method, whereas for
optimized and ESP settings, as well as all GROMACS settings, the AD method
is used (specifically by taking the analytic derivatives
of the piecewise polynomial approximants of $\phi$).

The approximation error
of the above five-step method for evaluating the truncated spectral sum
(the last formula in Eq.~\eqref{spectralpart}) is due to aliasing
induced by sampling on the $h$ grid~\autocite{hockney1988}.
It can be shown that the dominant aliasing error in the spectral part of the potential is of the order
\be
E_A= 6 \sum_{\ki\in I, \ki\neq\mbf{0}}
\left|\frac{\hat\chi(2\pi\rc |\ki|/L)}{|\ki|^2 \,\hat\phi(2\pi\ki/L)}
\,\hat\phi\left(\frac{2\pi}{L}[\ki + (n_f,0,0)]\right)\right|~,
\label{errorformula}
\ee
where the factor of 6 accounts for the six nearest-neighbor images in the reciprocal cubical lattice (displaced by $\pm 2\pi/h$ along each Cartesian axis). 
The total error in the potential is the sum of the truncation error in \eqref{truncationerror} and the aliasing error in \eqref{errorformula}. Detailed analysis is somewhat involved and will be reported at a later date.  

We sketch general rules for parameter selection as follows.
For the precision one generally sets $\eps = \Delta$, the target relative $\ell_2$ force error
\eqref{Deltadef}.
The shape parameter $\alpha$ is first determined by setting $\chi_\alpha(1)=\eps$
to control local truncation error.
For a given  cutoff radius $\rc$, a minimum $n_f$ is determined by
$\hat{\chi}(\pi \rc n_f/L) \le \eps$,
to ensure small truncation error in the spectral sum \eqref{spectralpart}.
The SI kernel is then chosen, and $n_f$ possibly increased, until
\eqref{errorformula} is of order $\eps$.
Finally, $\rc$ may be adjusted to balance the computational costs of the local part $u_i^l$ \eqref{localpart} and the spectral part $u_i^s$ \eqref{spectralpart}. The evaluation of $u_i^l$ scales as $O(Ns)$, where the average number of particles within the cut-off radius $s \propto \rc^3$. In contrast, the spectral evaluation involves an $O(NP^3)$ cost for spreading and interpolation, while the FFT cost scales as $O(N_F \log N_F)$ with $N_F \propto \rc^{-3}$. While the specific prefactors for these terms depend on the underlying computer architecture, the inverse dependence on $\rc^3$ ensures that $\rc$ can always be tuned to optimize the performance of the fast Ewald summation.

\subsection*{Classical fast Ewald methods}

In classical fast Ewald methods, following~\autocite{hockney1988},
the splitting kernel $\chi$ in \eqref{kernelsplitting} and the SI 
kernel $\phi$ for spreading/interpolation are chosen via
\be
\chi(r) = \frac{2}{\sqrt{\pi}}e^{-\alpha r^2},
\qquad \hat{\phi}(\kk) = \prod_{i=1}^3\hat{\phi}(\xi_i), \quad
\hat{\phi}(\xi_i)=\left(\frac{\sin(\xi_ih/2)}{\xi_ih/2}\right)^P,
\label{cem}
\ee
in both LAMMPS (PME) and GROMACS (PPPM).
This corresponds to the spatial SI kernel $\phi$ being a B-spline of width $P$ grid points per dimension~\autocite{darden1993jcp,ballenegger2012}.
The shape parameter $\alpha = \al_{\text{G}}$
for the Gaussian kernel,
and the minimal number of Fourier modes $n_f$ in each dimension are given by
\be
\alpha_{\text{G}} = \log(1/\eps), \qquad 
n_f =2\left \lceil \log(1/\eps)\frac{L}{\pi \rc}\right \rceil,
\label{cemparameters}
\ee
via the real- and Fourier-space truncation criteria, respectively.
However, due to the rather low default value $P=5$ in
both LAMMPS and GROMACS,
limitations imposed by the aliasing error \eqref{errorformula}
often require $n_f$ to be pushed significantly greater than the value listed in     \eqref{cemparameters}.
For MD simulations, this increase in $n_f$
leads to much larger FFTs and hence
communication costs~\autocite{ayala2021scalability}.

\subsection*{Optimizing the evaluation of the spectral part using PSWFs}

We now show that using PSWFs (the ESP method) offers a significant reduction
in the number of Fourier modes, leading to both  greater efficiency and better parallel scalability. 
For this, we denote the first prolate spheroidal wave function (PSWF) of order zero by $\psi_0^c$, where $c$ is its shape parameter,
and recall that it is defined~\autocite{slepian1961bstj}
as the eigenfunction with the largest eigenvalue
of the compact integral operator $\mathscr{F}_c:L^{2}[-1,1]\rightarrow L^2[-1,1]$,
defined via the formula
\be
\mathscr{F}_c[\varphi](x)=\int_{-1}^{1}\varphi(t)e^{icxt}dt.
\ee
We use a PSWF for both the splitting kernel and its Cartesian
products for the SI kernel:
\be
\chi(x) = \frac{1}{C_0}\psi_0^c(x),
\qquad \phi(\x) = \prod_{i=1}^3\phi(x_i), \quad
\phi(x_i)=\psi_0^{c_1}(2x_i/(Ph)),
\label{pswf}
\ee
where the normalization factor is $C_0=\int_{-1}^1 \psi_0^c(x)dx$.
As before, $P$ sets the SI kernel width in grid points in each dimension.

While the advantages of replacing Gaussians with PSWFs have been understood
in signal processing and imaging since their initial discovery,
many decades passed before they began to play a role in scientific computing.
This was due in part to the lack of simple closed-form expressions for PSWFs and their historically limited adoption in the MD community. That is no longer a barrier, as the evaluation of these functions is now efficiently performed via piecewise polynomial approximations. By using precomputed coefficients---consistent with the handling of tabulated potentials and spreading kernels in GROMACS and LAMMPS---these functions are evaluated with the same high efficiency as traditional methods.
The advantage of replacing the Gaussian with a PSWF as an {\em Ewald splitting kernel}
was noted in the DMK framework\autocite{jiang2025dmk} but, as far as we know,
has not previously been coupled with FFT-based Ewald summation.

\begin{figure}[!htbp]
\centering
\includegraphics[height=60mm]{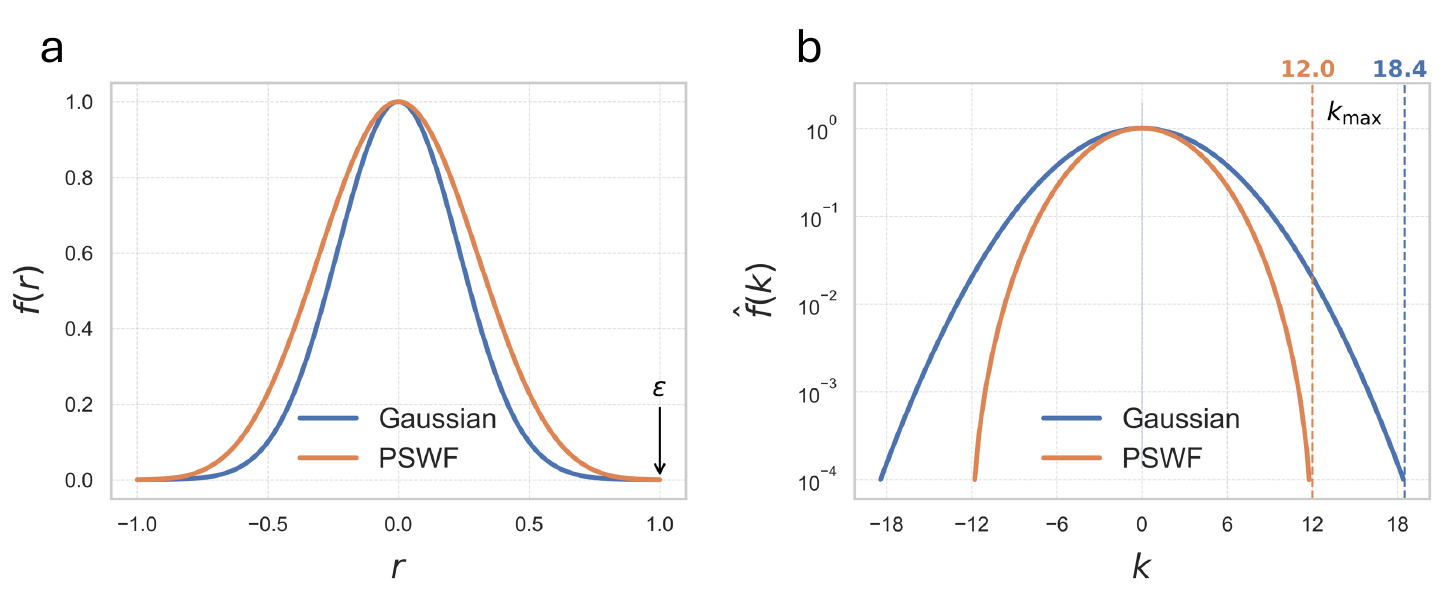}
\caption{\sf Comparison of Gaussian and prolate spheroidal wave function (PSWF) kernels used in Ewald splitting. {\bf a}, The Gaussian kernel $e^{-\alpha x^2}$ with $\alpha=9.2$, and
  the PSWF kernel $\psi_0^c(x)$ with $c=12$ plotted in real space. 
{\bf b}, Their Fourier transforms truncated
 to $\varepsilon=10^{-4}$. The dashed vertical lines mark the Fourier cutoffs used for the Gaussian and PSWF kernels, respectively.
The number of Fourier modes required by the PSWF is a factor of 
$(18.4/12.0)^3 \approx 3.6$ times smaller than when using the Gaussian for
Ewald splitting. Source data are provided as a Source Data file.}
\label{pswfgaussian}
\end{figure}

Turning to parameter choices,
we set $c$ such that $\psi_0^c(1)=\eps$,
and note that in the small-$\eps$ limit this tends to $c \approx \log(1/\eps)$.
The PSWF has the property that the
$\eps$-support of its Fourier transform is $[-c,c]$,
and recalling \eqref{hS} and the Fourier truncation at wavenumber $\pi n_f/L$,
we get
\be
c \approx \log(1/\eps), \qquad 
n_f =\left \lceil \frac{cL}{\pi \rc} \right \rceil, \qquad \mbox{ as } \eps\rightarrow 0,
\label{pswfparameters}
\ee
as the optimal parameters.
When compared with \eqref{cemparameters},
the number of necessary Fourier modes is reduced by a factor of $2$
in each dimension, in this high precision limit,
leading to a factor
of $8$ reduction in the total number of Fourier modes $N_f$ in the MD context.
At lower precisions the reduction is less dramatic,
but at 4 digits of accuracy, $N_f$ is still reduced by a factor of $3.6$.
In Fig.~\ref{pswfgaussian}, we contrast the Gaussian kernel and the PSWF,
chosen for $\eps=10^{-4}$, in both the physical and Fourier domains.
In the Fourier domain it is obvious that the PSWF may be truncated (to
accuracy $\eps$) at a significantly smaller frequency.

By a similar reasoning \autocite{finufft},
PSWFs are also the optimal window kernel for spreading and interpolation, minimizing
the number $P^3$ of grid points that couple 
to each particle. In fact, \eqref{errorformula}
suggests that if the SI kernel is the same physical width as the
Ewald split kernel, i.e., $Ph/2 = \rc$, the aliasing error is also $\bigO(\eps)$.
In practice, the value 
$c_1$ for the SI kernel $\phi$ in \eqref{pswf}
is better chosen slightly larger than the value $c$ in the splitting kernel
$\chi$.
This follows after extensive
numerical optimization at various precisions,
guided by the error formula \eqref{errorformula}.
Table~\ref{optimalvalues} lists recommended parameter values for the proposed ESP method.
Note that for our results in Table~\ref{tbl:expts} we sometimes found it
advantageous (depending on core count and other factors)
to vary the order $P$ from $P_{\rm esp}$ by at most 1.

\captionsetup{skip=20pt}
\renewcommand{\arraystretch}{1.35}
\begin{table}[!ht]
\caption{\sf Optimal parameter values for our prolate spheroidal wave function (PSWF)-based method (ESP). Here, $\eps$ represents the prescribed error tolerance of the electrostatic force. $\alpha_{\text{G}}=\log(1/\eps)$ is the Gaussian shape parameter in the Ewald splitting.
$c_{\rm pswf}$ is the value of the shape parameter $c$ of the PSWF $\psi_0^c$ in the ESP. $P_{\rm esp}$ is the optimal number of grid points per direction coupled to each particle in the ESP method. The last row, the quantity $R_{\rm default}$ denotes the ratio of the total number of Fourier grid points required in the PME method using default settings (with $n_f$ increased to achieve the stated error), to that in the ESP method
(with the same error).
Both LAMMPS and GROMACS set $P=5$ by default when running in MPI parallel; the achieved reduction in FFT grid size is dramatic, especially at higher precision.
}
\centering
\begin{tabular}{|c||c|c|c|c|c|}
\hline
$\eps$ &  $1\times 10^{-3}$ & $5\times 10^{-4}$  & $1\times 10^{-4}$ & $5\times 10^{-5}$ & $1\times 10^{-5}$\\
\hhline{|=||=|=|=|=|=|}
$\alpha_{\rm G}$ & $6.9078$ & $7.6009$ & $9.2103$ & $9.9035$ & $11.5129$\\
\hline
$c_{\rm pswf}$ & $9.5392$ & $10.290$ & $12.024$ & $12.762$ & $14.471$\\
\hline $P_{\rm esp}$ & $5$ & $5$ & $6$ & $7$ & $8$\\
\hhline{|=||=|=|=|=|=|}
$\mathbf{R_{\mathbf{ default}}}$ & $\mathbf{5.78}$ & $\mathbf{7.71}$ & $\mathbf{14.32}$ & $\mathbf{22.42}$ & $\mathbf{54.57}$\\
\hline
\end{tabular}
\label{optimalvalues}
\end{table}

\subsection*{Experimental design}

The pure-Coulomb performance tests on water, membrane, ribosome, truncated-octahedral CONA, and rhombic-dodecahedral RNA systems
(part 1 of Results) were set up as follows.
The water systems were built by isotropic periodic replication of a cubic box of side length $0.31$nm, containing a single water molecule. The membrane system uses a cuboid box of size $10.79$nm$\times 10.18$nm$\times9.55$nm. The ribosome system uses a cuboid box of size $31.25$nm$\times 31.25$nm$\times22.10$nm. For the CONA system, a triclinic cell equivalent to a truncated octahedron is used~\autocite{bekker1997unification}, with three equal edge lengths ($17.37$nm each) and inter-edge angles of $70.53^\circ$, $109.47^\circ$, and $70.53^\circ$. For the RNA system, a triclinic cell equivalent to a rhombic dodecahedron is used, with three equal edge lengths ($25.85$nm each) and inter-edge angles of $60^\circ$, $60^\circ$, and $90^\circ$. Periodic boundary conditions are imposed on all systems. For the membrane and ribosome systems, the input files and force field setups were taken from a recent benchmark~\autocite{kutzner2025scaling}, in which a cutoff of $1.0$ nm and Fourier spacings of $0.12$~nm (membrane) and $0.135$~nm (ribosome) were used for PME. For the water systems, we also use a cutoff of $1.0$~nm for consistent comparison, with Fourier spacings of $0.13~$nm and $0.083~$nm for $\Delta=10^{-3}$ and $10^{-4}$, respectively. For the triclinic systems, the input files and force field setups were taken from two application papers~\autocite{musleh2024analysis,posani2025ensemble}, in which cutoffs of $1.2~$nm (CONA) and $1.0~$nm (RNA) and Fourier spacings of $0.12$~nm (CONA) and $0.09$~nm (RNA) were used for PME. In modern MD software, LJ interactions are typically truncated using the same cutoff as the local part of the Coulomb interactions, allowing them to be treated efficiently together. To generate Fig.~\ref{fig:time} and Supplementary Fig.~7, all LJ interactions were disabled so that only Coulomb contributions were measured. For ESP only, the cutoff radius for the local part, $r_c$, was adjusted to optimize performance,
  resulting in the parameters in row 1--4 of Table~\ref{tbl:expts}.

For bulk water MD simulations (part 2 of Results),
we use the classical three-site SPC/E model, with fixed charges and Lennard--Jones (LJ) parameters assigned to each of the three atoms. All water MD simulations reported here were performed in the NVT ensemble. In all our MD simulations of water systems,
the short-range Coulomb and LJ potentials are truncated at $0.9$ nm, consistent with recommended settings to ensure sufficient accuracy of the LJ interactions. Four cubic simulation boxes of increasing size---$32.90$ nm (3,597,693 atoms), $34.82$ nm (4,214,784 atoms), $102$ nm (106,238,712 atoms), and $124.35$ nm (192,000,000 atoms)---were used to generate Figs.~\ref{fig:LAMMPS_time}(a--d) and Supplementary Fig.~2(a--d), respectively, all using periodic boundary conditions.
The LAMMPS water systems were constructed from $11^3$ and $34^3$ periodic copies of
a periodic box containing $901$ water molecules that had been equilibrated for 100 ns in the NVT ensemble~\autocite{frenkel2001} at $T = 298$ K.
The GROMACS water systems were constructed from $112^3$ and $400^3$ periodic copies of a single
water molecule from the standard GROMACS database, yielding a uniform initial density distribution.
These were followed by production runs using PPPM-, PME-, or ESP-based MD for data collection.
Temperature control and time integration were performed using the Nos\'e--Hoover chain thermostat\autocite{martyna1992nose}, with a time step of 2 fs and a coupling time constant $\tau = 100$ fs. All bond constraints were applied using the SHAKE algorithm~\autocite{krautler2001fast} in LAMMPS and the LINCS algorithm~\autocite{hess1997lincs} in GROMACS. \tcb{The parameter sets listed in Table~\ref{optimalvalues} correspond to the optimally tuned ESP parameters. In our benchmarks, the results labeled as ``default'' were obtained with spreading/interpolation kernel width $P=4$ for single-core tests and $P=5$ for multi-core tests (the default LAMMPS/GROMACS settings)}, respectively, then increasing the FFT grid size until the desired force accuracy was achieved.

For long-time MD simulations (part 3 of Results),
the lysozyme system uses the CHARMM27 force field for proteins~\autocite{foloppe2000all}, the TIP3P water model~\autocite{price2004modified}, and a cubic simulation box of size $22.02$ nm with periodic boundary conditions. The system comprises 1,036,152 atoms, including $27$ lysozyme proteins (1,960 atoms each), 327,672 water molecules, and $0.1$ M NaCl to reflect physiological ionic strength. The system is equilibrated for $500$ ns in the NPT ensemble at $298$ K and $1$ bar using the C-rescale barostat~\autocite{bernetti2020pressure}, followed by a $100$ ns production run in the NVT ensemble using the Nos\'e--Hoover chain thermostat and a time step of $2$ fs. Simulations are performed using GROMACS with either PME or ESP-based electrostatics.

The Li-ion aqueous electrolyte system uses the OPLS-AA force field for Li$^+$ ions~\autocite{jorgensen1996development}, the TIP3P model for water~\autocite{price2004modified}, and a specialized force field for TFSI$^-$ anions~\autocite{canongia2004molecular}. The system contains 1,011,392 atoms, including 20,480 Li$^+$, 20,480 TFSI$^-$, and 227,904 water molecules, and is simulated in a periodic cubic box of size $22.92$ nm using GROMACS. Equilibration is performed in the NPT ensemble at $298$ K and 1 bar for $200$ ns, followed by $100$ ns of production dynamics in the NVT ensemble using the Nos\'e--Hoover chain thermostat. The cutoff for both the short-range part of Coulomb, and LJ interactions, is set to $0.9$ nm. All parameter settings, such as the number of spreading/interpolation points and the real-space grid size, are consistent with those used in above benchmarks.

The relative RMS error of a computed set of electrostatic forces $\{\tilde{\F}_i\}_{i=1}^N$
is defined as usual as
\begin{equation}
  \Delta := \left( \frac{ \sum_{i=1}^N \| \tilde{\F}_i - \F_i \|_2^2 } { \sum_{i=1}^N \| \F_i \|_2^2 }
  \right)^{1/2},
\label{Deltadef}
\end{equation}
where $\F_i$ is an accurate reference calculation. In all measurements of $\Delta$ reported,
we computed the reference with estimated 8--9 digit accuracy by hand-tuning an increased spreading width, $P$,
and FFT grid size per dimension, $n_f$, in the same code.

It should be noted that the Gaussian shape parameter $\alpha_{\text{G}}$ used in this work, as listed in Table~\ref{optimalvalues}, is dimensionless and differs in definition from the corresponding ``G vector'' ($\alpha_{\text{PPPM}}$) in LAMMPS and the ``Gaussian width'' ($\alpha_{\text{GRO}}$) in GROMACS, both of which appear in the respective log files. The conversion between these definitions is straightforward: $\alpha_{\text{G}} = (\alpha_{\text{PPPM}}r_c)^2 =(r_c/\alpha_{\text{GRO}})^2$ which allows one to relate the dimensionless Gaussian shape parameter to the splitting factors used in the respective software, where $r_c$ is the real-space cutoff. Through our numerical tests, we find that the parameter selection scheme proposed in this work accurately meets the prescribed relative error tolerance (as in Table~\ref{optimalvalues}), whereas the default strategies implemented in LAMMPS and GROMACS often underestimate the actual relative error.
After consultation with the software developers, we note that the default settings in these packages are considered optimal for PME and PPPM in most practical simulations. Therefore, we adopt the default settings for comparison.

\subsection*{Methods for calculating physical properties}

The RMSF of particle $i$ is defined as the standard deviation of its position after least-squares fitting to a reference structure. The RMSF of a residue is computed as the average RMSF of all atoms within that residue.
The SASA of the protein molecule is calculated by rolling a sphere with a radius of the solvent probe of $1.4$ $\mathring{\text{A}}$ over the surface of the protein. RG denotes the radius of gyration  of the protein backbone as a function of time, where the atoms are explicitly mass weighted. The RMSD of a certain structure at time $t$ to a reference structure is calculated by least-squares fitting the structure to the reference structure and subsequently calculating the RMSD as 
\begin{equation}
f_{\text{RMSD}}(t)=\sqrt{\frac{1}{M}\sum_{i=1}^{N}m_i\left|\mathbf{r}_i(t)-\mathbf{r}_i(0)\right|^2},
\end{equation}
where $M=\sum_{i=1}^{N}m_i$ and $\mathbf{r}_i(t)$ is the position of atom $i$ at time $t$. The isochoric heat capacity ($C_{\text{Heat}}$) is calculated by the equation $C_{\text{Heat}}=\sigma_{\text{E}}^2/(2k_{\text{B}}T)$, where $\sigma_{\text{E}}$ is the standard deviation of the total energy in the system in the NVT ensemble, $k_{\text{B}}$ is the Boltzmann constant, and $T$ is the temperature. The dielectric constant is calculated by the equation
\begin{equation}
\varepsilon_{r}=1+\frac{4\pi}{3k_{\text{B}}TV}\left(\langle D^2\rangle-\langle D \rangle^2\right),
\end{equation}
where $V$ is the volume and $D$ is the total dipole moment of the system.

Native analysis tools in GROMACS permit the evaluation of RMSF, SASA, RMSD, RG, inter-residue distance distributions, time-averaged counts of residues by secondary structure type, the Gibbs energy difference, the Ramachandran distributions of residues, electrostatic potential and mass density profiles of proteins and lipids along the $z$-axis, RDFs, probability distributions of H--D--A angles and HB distances, dielectric constant, diffusion constant, and heat capacity. The charge structure factor is defined as the ensemble average in reciprocal space~\autocite{hu2022symmetry}: 
\begin{equation} 
S(\mathbf{k}) = \frac{1}{V} \sum_{i,j=1}^{N} q_i q_j \left\langle e^{i \mathbf{k} \cdot \mathbf{r}_{ij}} \right\rangle, 
\end{equation} 
where $\langle \cdot \rangle$ denotes an ensemble average. This quantity is computed using our own Python-based analysis tool. 

\subsection*{Hardware and software}
The computations in this paper were run on the Flatiron Institute Rusty cluster supported by the Scientific Computing Core at the Flatiron Institute. Each CPU node contains two AMD EPYC 9474F 48-core processors and 1.5 TB memory. The computing networks are connected using HDR-200 Infiniband (200 Gb/s). The GNU Compiler GCC/11.4.0 with OpenMPI/4.0.7 is used as the compiler for GROMACS, and FFTW/3.3.10 is used as the underlying FFT library.
The software version numbers were: LAMMPS version 19 Nov 2024, and GROMACS version 2025.1.

\section*{Data Availability}
\begin{sloppypar}
The datasets and input files for testing the PME, PPPM, and ESP methods---including SPC/E bulk water, lysozyme protein, and Li-ion aqueous electrolyte systems---have been deposited in the Zenodo database at DOI: 10.5281/zenodo.19559644. These materials are also available in the actively maintained GitHub repository
\url{https://github.com/LiangJiuyang/Ewald-Splitting-with-Prolates},
in the directory \verb+gromacs_bench+. Source data are provided with this paper. All data generated in this study are provided in the Zenodo link and/or Source Data files.
\end{sloppypar}

\section*{Code Availability}
\begin{sloppypar}
All original code has been deposited on Zenodo at DOI: 10.5281/zenodo.19547608~\autocite{liang2025_esp_code}. The active development repositories are available on GitHub:
LAMMPS-ESP at \url{https://github.com/LiangJiuyang/Ewald-Splitting-with-Prolates}
and GROMACS-ESP at
\url{https://github.com/lu1and10/Ewald-Splitting-with-Prolates}.
\end{sloppypar}

\printbibliography

\section*{Acknowledgments}
The authors are grateful for discussions with Pilar Cossio and Berk Hess.
They would like to thank the Scientific Computing Core at the Flatiron Institute for their support and for providing essential computational resources.
The work of J.L. is partially supported by the National Natural Science Foundation of China (Grant No. 12401570) and the China Postdoctoral Science Foundation (Grant No. 2024M751948).
The Flatiron Institute is a division of the Simons Foundation.

\section*{Author Contributions Statement}
S.J., A.B., and L.G. designed the project and contributed to overall project guidance. S.J. and J.L. designed the experiments. J.L. and L.L. implemented the ESP method within the LAMMPS and GROMACS framework and performed the simulations. All authors discussed the results and wrote the paper. All authors approved the final version.

\section*{Competing Interests Statement}
The authors declare no competing interests.

\newpage

\clearpage
\renewcommand{\thepage}{S\arabic{page}}
\setcounter{page}{1}

\subsection*{Supplementary Information}
\renewcommand{\figurename}{Supplementary Figure}
\setcounter{figure}{0} 

\begin{figure}[!htbp]
\centering
\includegraphics[width=0.98\linewidth]{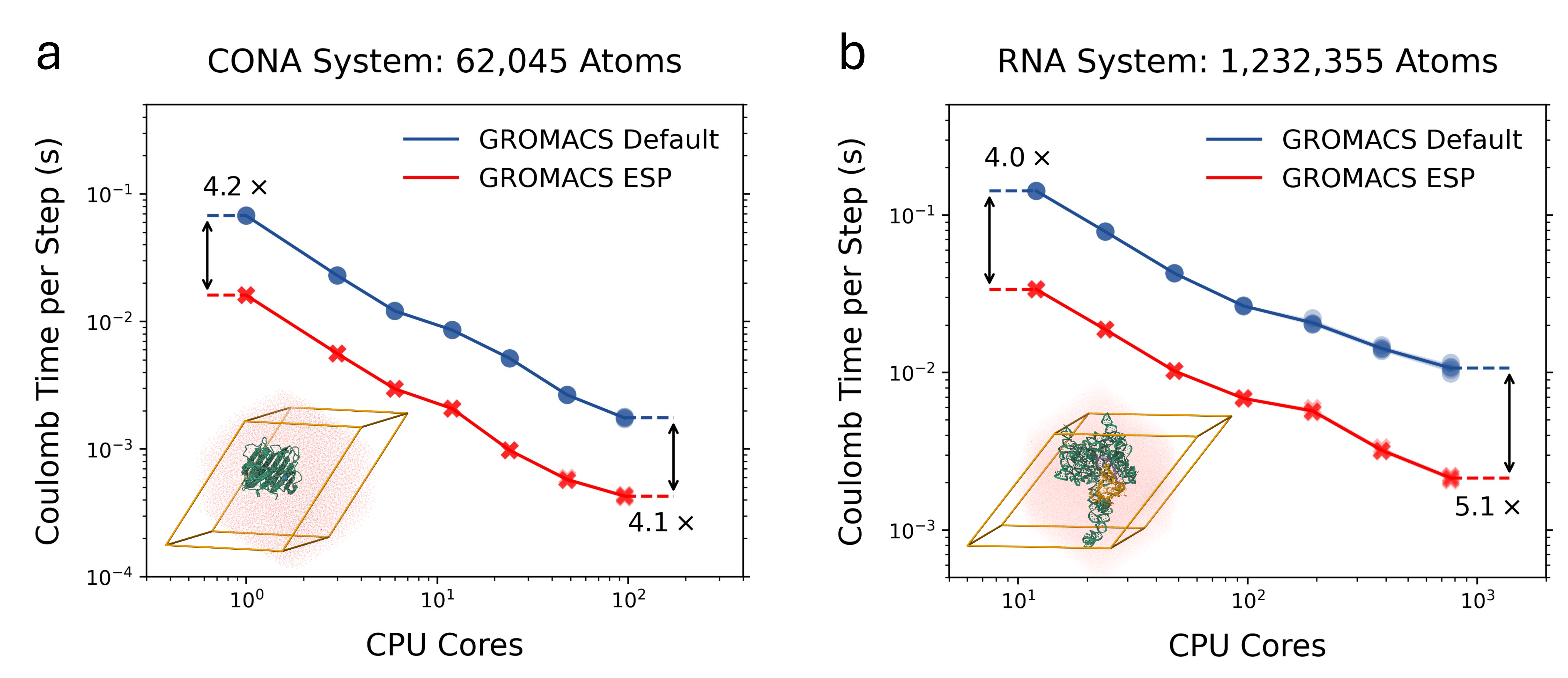}
\caption{\sf Performance comparison of ESP with PME for truncated-octahedral CONA (62,045 atoms) and rhombic-dodecahedral RNA (1,232,355 atoms) systems. Insets show snapshots overlaid with the equivalent triclinic unit cells used in the simulations (yellow), with water molecules shown in light red and the solute colored by chain index. Particles shown outside the displayed triclinic cells are wrapped back into the unit cells under periodic boundary conditions during simulations. The average Coulomb computation time per step is plotted as a function of the number of CPU cores, with each point averaged over five runs of $30$ minutes. The PME reference data in Panels {\bf a} and {\bf b} follow setups reported in~\cite{musleh2024analysis} and \cite{posani2025ensemble}, with estimated error tolerances of $4\times 10^{-4}$ and $2\times 10^{-4}$ (validated in Supplementary Fig.~7), respectively. Source data are provided as a Source Data file.}
\label{fig:CoulombOnlyPerformanceTriclinic}
\end{figure}

\begin{figure}[!htbp]
\centering
\includegraphics[width=0.98\linewidth]{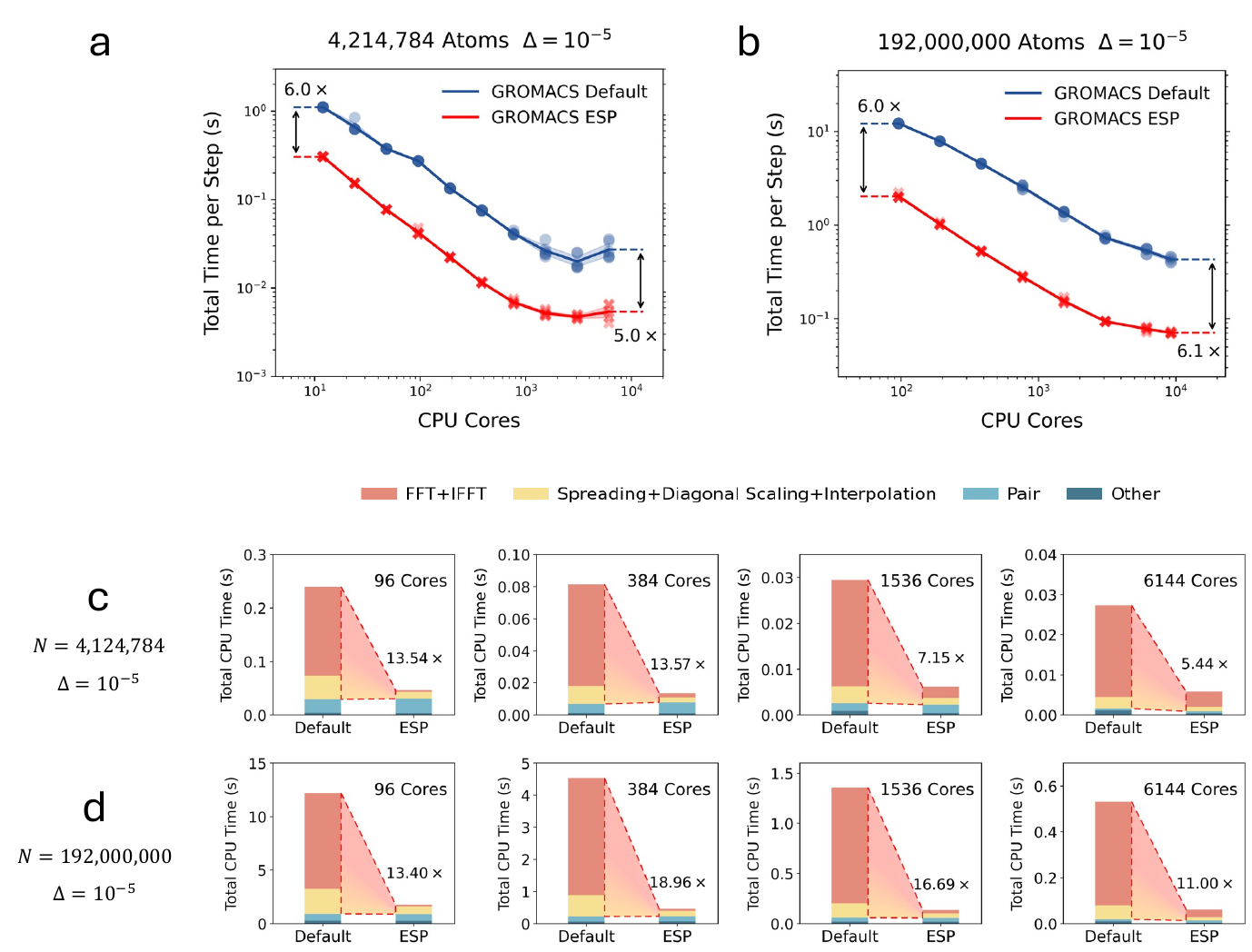}
\caption{\sf Performance comparison of ESP with PME for bulk water ({\bf a}, 4,214,784 atoms; and {\bf b}, 192,000,000 atoms)
systems. The total simulation time per step, averaged over five runs of $30$ minutes each, is shown for these systems as a function of the number of CPU cores. Data were generated using ESP implemented within GROMACS, 
compared against the native PME option, with an error tolerance of $\Delta = 10^{-5}$. 
Blue circles show PME with default parameters, and red '\texttimes{}' marks show ESP. Light-colored markers show the results of five repeated runs. Solid lines indicate the mean across the five runs, and shaded bands indicate the $95\%$ confidence intervals of the mean. In panels {\bf c} and {\bf d}, the simulation time per step is broken down into four components: FFT and IFFT operations (red); spreading, diagonal scaling, and interpolation (yellow); short-range pairwise interactions (light blue); and all other simulation tasks (dark blue). Panels {\bf c}-{\bf d} correspond to the systems shown in {\bf{a}}-{\bf b}, respectively. Speedups in the long-range interaction are annotated directly within each panel. Source data are provided as a Source Data file.}
\label{fig:GROMACSTime}
\end{figure}

\begin{figure}[!htbp]
\centering
\includegraphics[width=0.98\linewidth]{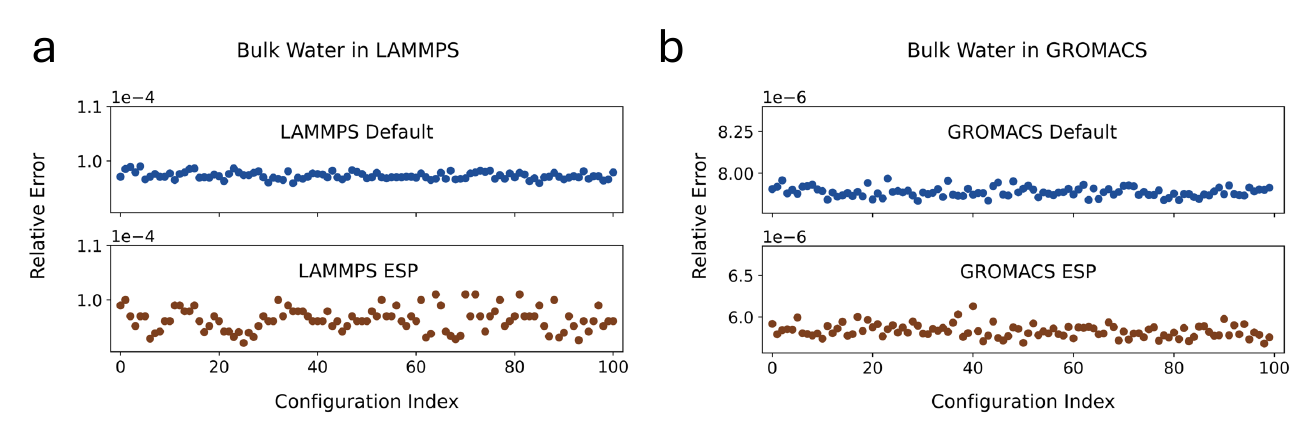}
\caption{Relative error across 100 bulk water configurations randomly sampled from the equilibration state is presented. The results are obtained using implementations in {\bf a}, LAMMPS with an error tolerance of $\Delta=10^{-4}$; and {\bf b}, GROMACS with an error tolerance of $\Delta=10^{-5}$. Comparisons are made between the proposed ESP variants of the codes based on PSWF, and the native LAMMPS/GROMACS codes under the default parameter setup. Source data are provided as a Source Data file.}
\label{fig:ErrorOverConfiguration}
\end{figure}

\begin{figure}[!htbp]
\centering
\includegraphics[width=0.98\linewidth]{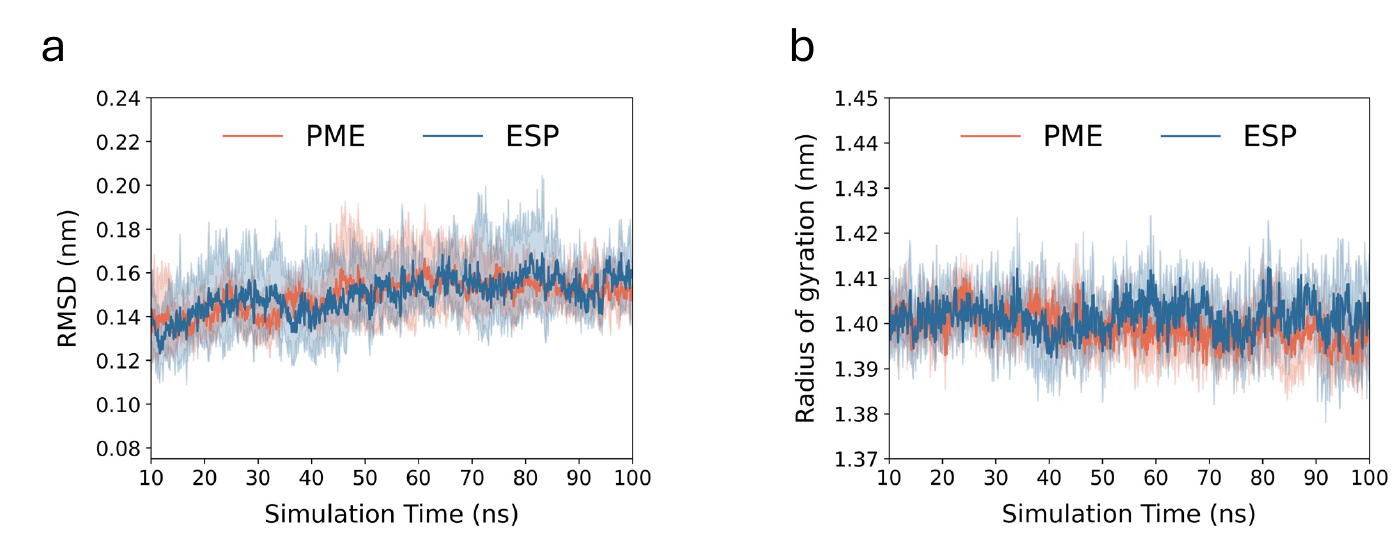}
\caption{Comparison of additional simulation results obtained using GROMACS with PME- and ESP-based MD on the same lysozyme protein system analyzed in Fig.~\ref{fig:protein}. {\bf a}, Root mean square displacement (RMSD) of the protein backbone. {\bf b}, Radius of gyration of the proteins. Data are shown for both the native PME-based GROMACS implementation and the ESP method over 100 ns of simulation. Light-colored shaded areas and colored markers in {\bf a}, {\bf b} indicate $95\%$ confidence intervals based on five independent runs. All results show good agreement between the two methods. Source data are provided as a Source Data file.}
\label{fig:ExtendStructualInformation}
\end{figure}

\begin{figure}[!htbp]
\centering
\includegraphics[width=0.98\linewidth]{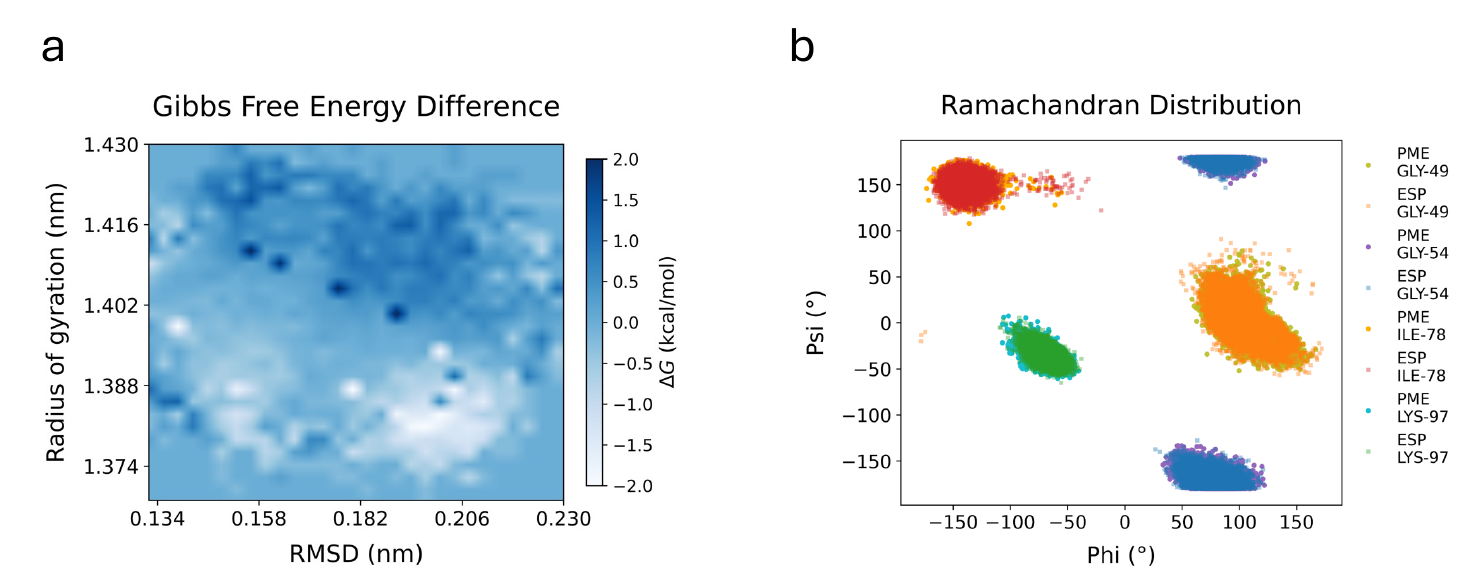}
\caption{Comparison of additional simulation results obtained using GROMACS with PME- and ESP-based MD on the same lysozyme protein system analyzed in Fig.~\ref{fig:protein}. The production phase lasts for $100$ ns. {\bf a}, Colormap representation of the Gibbs free energy difference between GROMACS and the ESP implementation in GROMACS, plotted against backbone RMSD ($x$-axis) and backbone radius of gyration ($y$-axis). The root mean square difference calculated across grid points with non-zero values is $0.6$ kcal/mol. {\bf b}, The Ramachandran distributions of glycine residues (GLY-49 and GLY-54), isoleucine residue (ILE-78), and lysine residue (LYS-97). All results show good agreement between the two methods. Source data are provided as a Source Data file.}
\label{fig:ExtendStructualInformation2}
\end{figure}

\begin{figure}[!ht]
\centering
\includegraphics[width=0.98\linewidth]{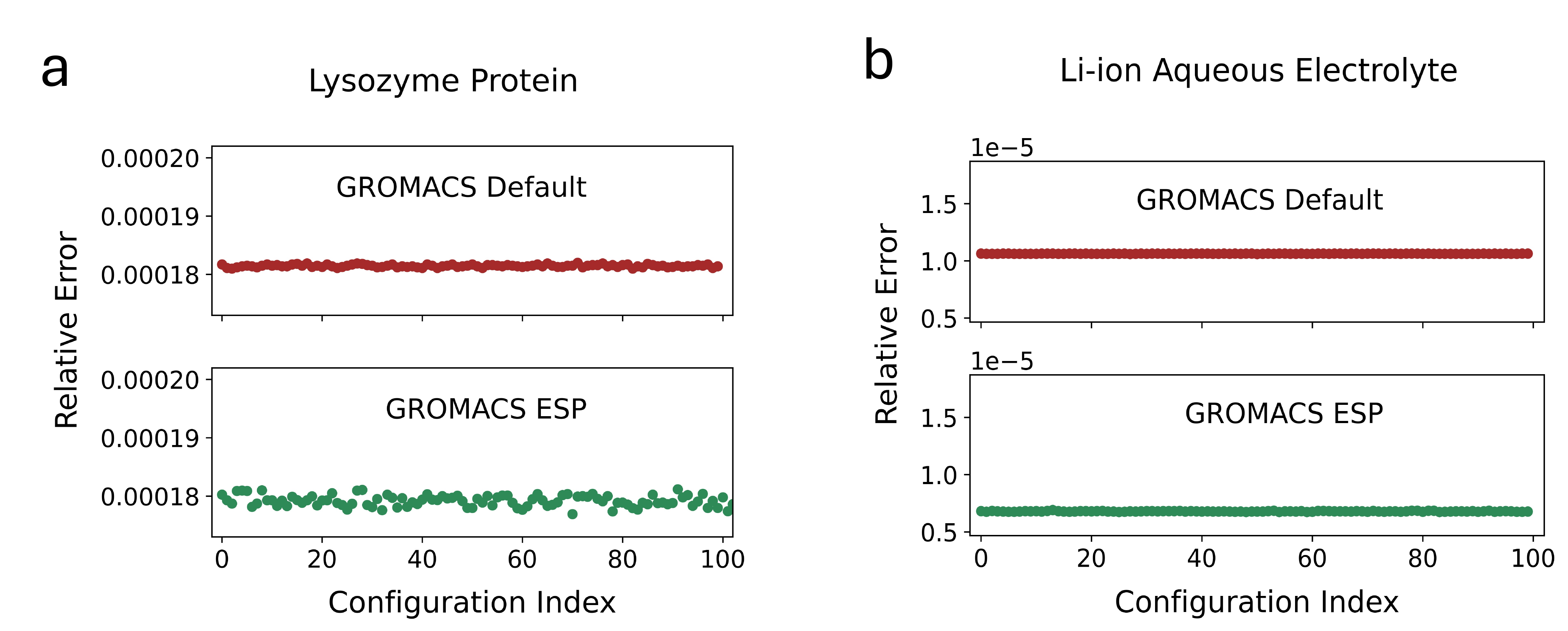}
\caption{Relative error across 100 configurations randomly sampled from the equilibration state. {\bf a}, lysozyme protein. {\bf b}, Li-ion aqueous electrolyte. The results are obtained using implementations in GROMACS with an error tolerance of $\Delta=2\times10^{-4}$ and $\Delta=10^{-5}$, respectively. Comparisons are made between ESP implemented in GROMACS, and the native GROMACS code under the default parameter setup. Source data are provided as a Source Data file.}
\label{fig::ErrorConf}
\end{figure}

\begin{figure}[!htbp]
\centering
\includegraphics[width=1\linewidth]{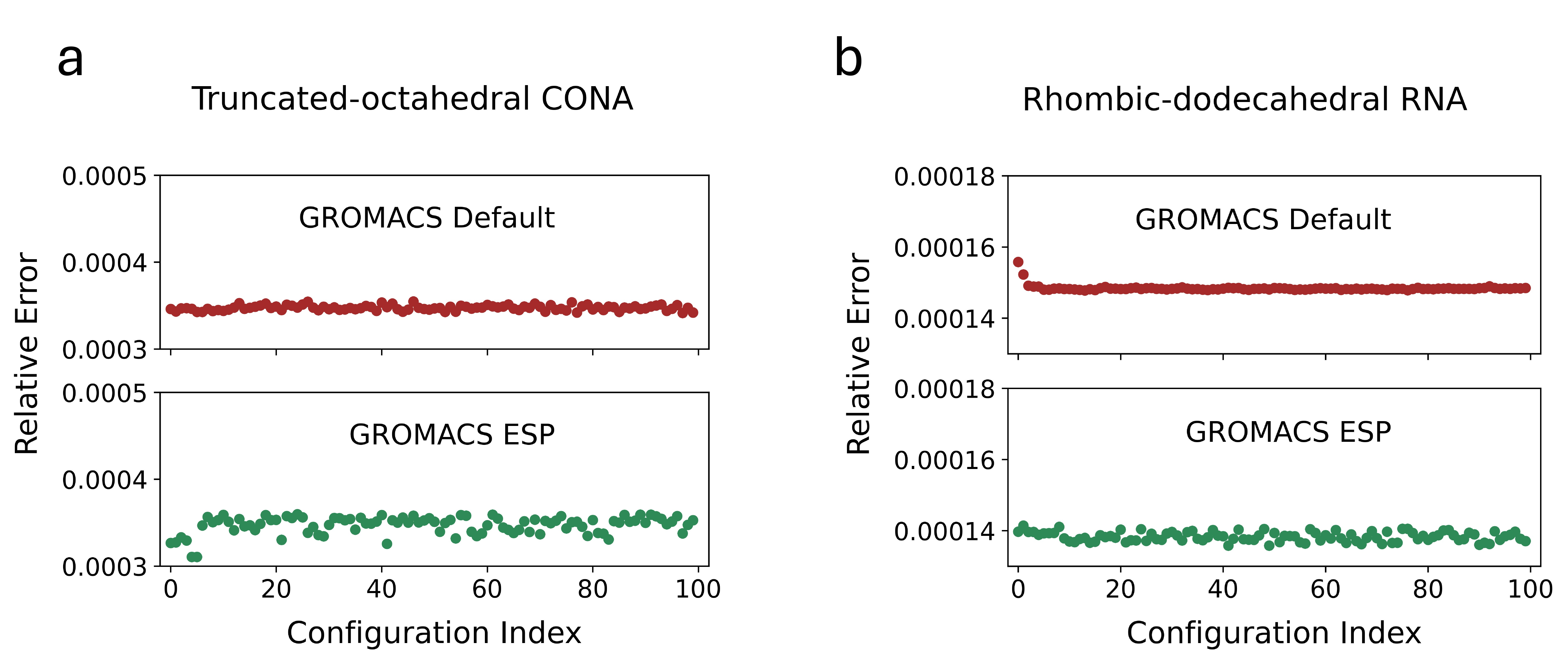}
\caption{Relative error across 100 configurations randomly sampled from the equilibration state is presented. {\bf a}, truncated-octahedral Concanavalin A (CONA). {\bf b}, rhombic-dodecahedral RNA. The results are obtained using implementations in GROMACS with an error tolerance of $\Delta=4\times 10^{-4}$ and $\Delta=2\times 10^{-4}$, respectively. Comparisons are made between the proposed ESP variants of the codes based on PSWF, and the native GROMACS code under the default parameter setup. Source data are provided as a Source Data file.}
\label{fig:ErrorOverTriclinicConfiguration}
\end{figure}


\end{document}